\documentclass[a4paper,12pt,nosumlimits]{amsart} 
\usepackage{amssymb}
\usepackage[english]{babel}
\usepackage{fullpage,xypic,enumerate,cancel}

\title{The Gysin Sequence for Quantum Lens Spaces}
\date{April 2014}

\author{Francesca Arici, Simon Brain, Giovanni Landi}
\address{SISSA, Via Bonomea 265, 34136 Trieste, Italy} 
\address{Matematica, 
Universit\`{a} di Trieste, Via A.~Valerio 12/1, 34127 Trieste, Italy}
\address{Matematica, 
Universit\`{a} di Trieste, Via A.~Valerio 12/1, 34127 Trieste, Italy 
and INFN, Sezione di Trieste, Trieste, Italy}
\email{farici@sissa.it, sbrain@units.it, landi@units.it}

\keywords{Noncommutative geometry, Gysin sequence, quantum lens spaces}
\subjclass[2010]{46L85; 58B34}

\DeclareFontFamily{OT1}{pzc}{}
\DeclareFontShape{OT1}{pzc}{m}{it}{<-> s * [1.20] pzcmi7t}{}
\DeclareMathAlphabet{\mathpzc}{OT1}{pzc}{m}{it}

\numberwithin{equation}{section}

\newtheorem{thm}{Theorem}[section]
\newtheorem{prop}[thm]{Proposition}

\newtheorem{defi}[thm]{Definition}
\newtheorem{exa}[thm]{Example}
\newtheorem{rem}[thm]{Remark}

\newtheorem*{prop*}{Proposition}

\newcommand{\U}{\mathrm{U}}
\newcommand{\A}{\mathcal{A}}
\newcommand{\B}{\mathcal{B}}
\newcommand{\HH}{\mathcal{H}}
\newcommand{\LL}{\mathcal{L}}
\newcommand{\N}{\mathbb{N}}
\newcommand{\Z}{\mathbb{Z}}
\newcommand{\R}{\mathbb{R}}
\newcommand{\C}{\mathbb{C}}
\newcommand{\CP}{\mathbb{C}\mathrm{P}}
\newcommand{\rL}{\mathrm{L}}
\newcommand{\rS}{\mathrm{S}}
\newcommand{\hp}[1]{\left<#1\right>}

\newcommand{\tr}{\mathrm{Tr}}

\newcommand{\coker}{\mathrm{coker}}
\newcommand{\diag}{\mathrm{diag}}
\newcommand{\beq}{\begin{equation}}
\newcommand{\eeq}{\end{equation}}
\newcommand{\nn}{\nonumber}
\newcommand{\Dom}{\mathfrak{Dom}}  
\newcommand{\ii}{\mathrm{i}}

\renewcommand{\tilde}{\widetilde}

\begin{document}

\begin{abstract}
We define quantum lens spaces as `direct sums of line bundles' and exhibit them as `total spaces' 
of certain principal bundles over quantum projective spaces. 
For each of these quantum lens spaces we construct an analogue of the classical Gysin sequence in K-theory. 
We use the sequence to compute the K-theory of the quantum lens spaces, in particular to give 
explicit geometric representatives of their K-theory classes. These representatives are interpreted 
as `line bundles' over quantum lens spaces
and generically define `torsion classes'. We work out explicit examples of these classes.

\end{abstract}

\maketitle

\tableofcontents

\linespread{1.1} 
\parskip 1ex

\vfill
{\footnotesize
SB was supported by INdAM, cofunded under the Marie Curie Actions of the European Commission (FP7-COFUND). GL was partially supported by the Italian Project ``Prin 2010-11 -- Operator Algebras, Noncommutative Geometry and Applications''.
}

\newpage
\section{Introduction}
This paper is devoted to the study of the noncommutative topology of quantum lens spaces via their K-theory. We construct an exact sequence --  a noncommutative analogue of the classical Gysin sequence -- which relates the K-theory of quantum lens spaces to the K-theory of quantum projective spaces. Our construction enables us not only to compute the K-theory of the quantum lens spaces in a novel way, but also to obtain geometric representatives of the K-theory classes, generically torsion ones, in terms of  `line bundles'.

Noncommutative (or quantum)  lens spaces have been the subject of increasing interest of late. They first appeared in  \cite{MT} in the context of what we would now call `theta-deformed' topology; they later surfaced in \cite{HS03} in the guise of graph $C^*$-algebras, with certain more recent special cases ({\em cf}. \cite{BF12,HRZ11}). The particular case of the quantum three-dimensional real projective space was studied in \cite{Po95} and \cite{La98}. 
Real spectral triples on three-dimensional noncommutative lens spaces have recently been  studied in \cite{S-V}.

Lens spaces arise in classical geometry as quotients of odd-dimensional spheres by an action of a finite cyclic group. 
In parallel with this, quantum lens spaces are usually introduced in terms of fixed point algebras for suitable actions of finite cyclic groups on function algebras over odd dimensional quantum spheres. Indeed, the key result of \cite{HS03} is the realization of the $C^*$-algebra of continuous functions on a quantum lens space as the Cuntz-Krieger algebra of a directed graph. From this,
and quite importantly, one deduces the K-theory of the algebra as the kernel and cokernel of a certain `incidence matrix' associated to the graph. This computation of the K-theory is, one has to say, very direct, although somewhat implicit and obtained via some rather 
complicated isomorphisms. 

In the present paper we task ourselves with finding a more elegant intuitive and geometric approach to the K-theory of quantum lens spaces. To this end, our starting point is the `algebraic' approach to the K-theory of quantum projective spaces presented in  \cite{DL10}. 
Center stage there is taken by (polynomial) bimodules $\LL_{N}$ of sections of noncommutative `line bundles' over the projective space; line bundles which determine the K-theory of the $C^*$-algebra $C(\CP^n_q)$. Out of this algebraic approach to K-theory there come several important advantages that we list in the remainder of this introduction, by way of summarizing some of the main results of the present paper. 

Given a pair of  positive integers $n,r$,  
the coordinate algebra $\A(\rL^{(n,r)}_q)$ of the \textit{quantum lens space} 
of dimension $2n+1$ (and index $r$) is defined to be 
\beq\tag{\ref{dals}}
\A(\rL^{(n,r)}_q) :=\bigoplus_{N\in\Z} \, \LL_{r N}\, .
\eeq
Then $\rL^{(n,r)}_q$ is the `total space' of a principal bundle over the quantum projective space $\CP^n_q$ with structure group $\widetilde{\U}(1) :=\U(1)/\Z_r$. This parallels the $\U(1)$ principal bundle over $\CP^n_q$ having total space the quantum sphere 
$\rS^{2n+1}_q$, the latter being obtained for $r=1$ in the previous decomposition:
\beq\tag{\ref{Lnbdl}}
\A(\rS^{2n+1}_q) :=\bigoplus_{N\in\Z} \, \LL_{N} \, .
\eeq

One is then able to show {\em a posteriori}  that the algebra $\A(\rL^{(n,r)}_q)$ is made of all elements of $\A(\rS^{2n+1}_q)$ which are invariant under a certain action of the cyclic group $\Z_r$.  With these principal bundles there comes a way to `pull-back'  line bundles from $\CP^n_q$ to $\rL^{(n,r)}_q$:  
\beq\tag{\ref{dia-pb}}
\xymatrix{ 
\widetilde{\LL}_N \ar@{.>}[d] & \ar[l]_{j_{*}} \LL_N \ar@{.>}[d] & \\  
 \A(\rL^{(n,r)_q}) \,\, & \,\, \ar[l]^j \A(\CP^n_q) \, .  &
}
\eeq
That is to say, the algebra inclusion $j : \A(\CP^n_q)\to \A(\rL^{(n,r)}_q)$ also induces a map
\beq\tag{\ref{Kpblb}}
j_* : K_0(C(\CP^n_q)) \to K_0(C(\rL^{(n,r)}_q)) \, .
\eeq
The marked difference between a line bundle $\LL_N$ over $\CP^n_q$ versus its pulled-back $\widetilde{\LL}_{N}$ to $\rL^{(n,r)}_q$ 
is that, while each $\LL_N$ is not free when $N\ne 0$, this need not be the case for  $\widetilde{\LL}_{N}$:
the pulled-back $\widetilde{\LL}_{-r}$ of $\LL_{-r}$ is tautologically free, that is to say it is trivial in the group 
$K_0(C(\rL^{(n,r)}_q))$. It follows that $(\widetilde{\LL}_{-N})^{\otimes r}\simeq\widetilde{\LL}_{-rN}$ also has trivial class for any $N\in\Z$ and thus such line bundles $\widetilde{\LL}_{-N}$ define \emph{torsion classes}; they generate the group $K_0(C(\rL^{(n,r)}_q))$. 

In addition, there is a multiplicative structure on the group $K_0(C(\CP^n_q))$: 
\begin{prop*}{\bf \ref{ku}.}
It holds that
$$
K_0(C(\CP^n_q)) \simeq \Z[ [\LL_{-1}] ] / (1-[\LL_{-1}])^{n+1}  \simeq \Z[u] / u^{n+1} \, 
$$
where $u = \chi([\LL_{-1}]) := 1 - [\LL_{-1}]$ is the Euler class of the line bundle $\LL_{-1}$.
\end{prop*}
\noindent
Out of this one is led to a map 
\beq\tag{\ref{Gsqls}}
\alpha : K_0(C(\CP^n_q)) \to K_0(C(\CP^n_q))
\eeq
where $\alpha$ is now multiplication by the Euler class $\chi(\LL_{-r}) := 1 - [\LL_{-r}]$ of $\LL_{-r}$. 
Central for us is the assembly of this map with the pull-back map \eqref{Kpblb} into an exact sequence
$$
\xymatrix{ 
0\to K_1(C(\rL^{(n,r)}_q)) \ar[r]^-{\textup{Ind}} & K_0(C(\CP^n_q)) \ar[r]^-{\alpha} 
& K_0(C(\CP^n_q)) \ar[r]^-{j_*} & K_0(C(\rL^{(n,r)}_q)) \ar[r] & 0,
}
$$
the \emph{Gysin sequence} for our quantum lens space $\rL^{(n,r)}_q$, with a suitable index map $\textup{Ind}$ explicitly described below. 
Having arrived at this sequence, one could easily be content simply by admiring its sheer elegance. 
It has, however, some very practical and doubtlessly important applications, which we present in the final sections.

Notably, there is the computation of the K-theory of the quantum lens spaces 
$\rL^{(n,r)}_q$. Owing to Prop.~\ref{ku},  the map $\alpha$ can be given as an $(n+1)\times (n+1)$ matrix with respect to the $\Z$-module basis 
$\{1, u, \dots, u^n \}$ of $K^0(C(\CP^n_q)) \simeq \Z^{n+1}$. This leads to the identifications
\beq\tag{\ref{qk1ls}}
K_1(C(\rL^{(n,r)}_q)) \simeq \ker(\alpha) , \qquad K_0(C(\rL^{(n,r)}_q)) \simeq \coker(\alpha) \, . 
\eeq
We stress that our construction is structurally different from the one in \cite{HS03}, the only point of contact being that  
the K-theory is obtained out of a matrix. First, our matrix is different from the incidence matrix of \cite{HS03}. Second, and more importantly, the structure of the map $\alpha$ and Prop.~\ref{ku} allow us to give \emph{geometric} generators of both the groups  
$K_1(C(\rL^{(n,r)}_q))$ and $K_0(C(\rL^{(n,r)}_q))$, for the latter in particular as (combinations) of pulled-back  line bundles from 
$\CP^n_q$ to $\rL^{(n,r)}_q$. All of this is described in full detail in \S\ref{ktqls}. 
 
Some of the dual constructions pertinent to the K-homology of the quantum lens spaces $\rL^{(n,r)}_q$, 
stemming from a sequence dual to the previous one, will be reported elsewhere.

\subsubsection*{Notation}
By a $*$-algebra we mean a complex associative unital involutive algebra. An unadorned tensor product is meant to be over $\C$.
As it is customary, a noncommutative $C^*$-algebra $A$ is thought of as being the algebra of continuous functions on an underlying `quantum' topological space, and we use the notation $K_{\bullet}(A)$ for the K-theory of this $C^*$-algebra, together with $K^{\bullet}(A)$ for its K-homology. 

\subsubsection*{Acknowledgments} 
We are grateful to Alan Carey, Francesco D'Andrea, Erik van Erp, Sasha Gorokhovsky, Jens Kaad, Max Karoubi, Ryszard Nest and 
Georges Skandalis for useful discussions. Adam Rennie deserves a special mention for making transparent one of our central theorems below. Finally, we thank an anonymous referee for some excellent comments which led to a much improved version of the paper. 
  
\section{The classical Gysin sequence}
In this section we simply follow Karoubi's book \cite{Ka78}. 
Let $\CP^n$ denote the complex projective space of $\C^{n+1}$ and let $V$ be a complex vector bundle over $\CP^n$ equipped with a Hermitian fibre metric. We write $B(V)$ for the `ball bundle' of $V$, the bundle over $\CP^n$ whose fibre $B(V)_x$ at the point $x\in\CP^n$ is the closed unit ball of the fibre $V_x$ of $V$. Similarly we write $S(V)$ for the `sphere bundle' of $V$, whose fibre $S(V)_x$ at $x\in\CP^n$ is the unit sphere of the fibre $V_x$. Then $B(V)-S(V)$ denotes the open ball bundle.

Since $S(V)$ is closed in $B(V)$, with $K^*(B(V),S(V))$ denoting the relative K-theory groups
one has a six term exact sequence in topological K-theory \cite[IV.1.13]{Ka78}:
\begin{equation}\label{six}
\xymatrix{ 
K^0(B(V),S(V))\ar[r] & K^0(B(V))\ar[r]  & K^0(S(V))\ar[d]^{\delta_{01}} \\  K^{1}(S(V))\ar[u]^{\delta_{10}} 
& K^{1}(B(V))\ar[l] & K^{1}(B(V),S(V)). \ar[l]
}
\end{equation}
Here the vertical arrows are the usual `connecting homomorphisms' , while
the horizontal arrows are induced by natural maps and will be described explicitly below \cite[II.3.21]{Ka78}.

Since $B(V)$ is compact, it follows that $K^*(B(V),S(V))\simeq K^\ast(B(V)-S(V))$. 
Moreover, the total space of the fibre bundle $B(V)-S(V)$ is homeomorphic to the total space of $V$. 
These facts, followed by the Thom isomorphism combined with Bott periodicity gives rise to isomorphisms of K-groups 
$$
K^*(B(V),S(V))\simeq K^\ast(B(V)-S(V))\simeq K^\ast(V)\simeq K^\ast(\CP^n).
$$ 

Finally, since the total space of $B(V)$ is homotopic to $\CP^n$ (via the inclusion of the latter into $B(V)$ determined by the zero section of $V$), one has isomorphisms of K-groups 
$$
K^\ast(B(V))\simeq K^\ast(\CP^n).
$$ 
Assembling all of this together and using the vanishing $K^1(\CP^n)=0$ (\emph{cf}. \cite[Cor.~IV.2.8]{Ka78}),
the sequence~\ref{six} transforms into the K-theoretic \emph{Gysin sequence} for the bundle $S(V)$:
\begin{equation}\label{five}
0 \longrightarrow K^1(S(V)) \xrightarrow{\,\,\delta_{10}\,} K^0(\CP^n) \xrightarrow{ \,\, \alpha \,\, } K^0(\CP^n) \xrightarrow{ \,\,  \pi^* \,} 
K^0(S(V)) \longrightarrow 0 \, .
\end{equation}
The `pull-back' homomorphism $\pi^*$ is induced by the bundle projection $\pi : S(V) \to \CP^n$. The homomorphism $\alpha$ is given by multiplication by the \emph{Euler class} $\chi(V)$ of the vector bundle $V$ that we shall momentarily describe more explicitly. 

Now let $L$ be the tautological line bundle over $\CP^n$, whose total space is $\C^{n+1}$ and whose fibre $L_x$ at $x\in \CP^n$ is the one-dimensional complex vector subspace of $\C^{n+1}$ which defines that point.  
Via the usual associated bundle construction, the bundle $L$ may be identified with the quotient of $\rS^{2n+1} \times\C$ by the equivalence relation 
$$
(x, t) \sim (\lambda x, \lambda^{-1} t),\qquad \lambda\in \rS^1 \subseteq \C.
$$ 
Similarly, its $r$-th tensor power $L^{\otimes r}$ may be identified with the quotient of 
$\rS^{2n+1} \times \C$ by the equivalence relation $(x, t) \sim (\lambda x, \lambda^{-r} t)$. Moreover, $L^{\otimes r}$ can be given the fibre metric defined by $\varphi\left( (x, t'),  (x, t) \right) = t' \bar{t}$. It follows that the sphere bundle $S(L^{\otimes r})$ can be identified with the `lens space' $\rL^{(n,r)} := \rS^{2n+1} / \Z_r$ (where the cyclic group $\Z_r$ of order $r$ acts upon the sphere $\rS^{2n+1}$ via the $r$-th roots of unity) by the map $(x, t) \mapsto \sqrt[r]{t} \cdot x $. 

Taking $V=L^{\otimes r}$ in the above sequence \eqref{five} one finds, just as in \cite[IV.1.14]{Ka78}, the K-theoretic \emph{Gysin sequence} for the lens space $\rL^{(n,r)}$:
\beq\label{cgs}
0 \longrightarrow K^1(\rL^{(n,r)}) \xrightarrow{\,\,\delta_{10}\,} K^0(\CP^n) \xrightarrow{ \,\, \alpha \,\, } K^0(\CP^n) \xrightarrow{ \,\,  \pi^* \,} 
K^0(\rL^{(n,r)}) \longrightarrow 0 \, .
\eeq
Here, since $L^{\otimes r}$ is a line bundle, its Euler class (giving the map $\alpha$) is given simply by
$$
\chi(L^{\otimes r}): = 1 - [L^{\otimes r}] \, .
$$  

\section{Quantum projective spaces}

We first describe the class of noncommutative projective spaces that we need. 
We recall both the algebras of coordinate and continuous functions on quantum projective space, together with the noncommutative `line bundles' which represent the K-theory.

\subsection{Functions on quantum projective spaces} In the following, without loss of generality, the real deformation parameter is restricted to the interval $0<q<1$. We recall from 
\cite{VS91} that the coordinate algebra  of the unit quantum sphere $\rS^{2n+1}_q$ is the $*$-algebra $\A(\rS^{2n+1}_q)$ generated by 
$2n+2$ elements 
$\{z_i,z_i^*\}_{i=0,\ldots,n}$ subject to the relations:
\begin{align}\label{eq:3rd}
z_iz_j &=q^{-1}z_jz_i && 0\leq i<j\leq n \;, \nn \\
z_i^*z_j &=qz_jz_i^*  &&  i\neq j \;, \nn \\
[z_n^*,z_n] &=0 \;, \qquad [z_i^*,z_i] =(1-q^2)\sum_{j=i+1}^n z_jz_j^* && i=0,\ldots,n-1 \;,  \nn \\
1 & = z_0z_0^*+z_1z_1^* +\ldots+z_nz_n^* \;.
\end{align}
The notation of \cite{VS91} is obtained by setting $q=e^{h/2}$,
while the relationship with the generators $x_i$ used in~\cite{HL04} is given by $x_i=z_{n+1-i}^*$ together with the replacement $q\to q^{-1}$.  
The quantum sphere $\rS^{2n+1}_q$ is a quantum homogeneous space of the quantum $\mathrm{SU}(n+1)$ group of \cite{Wor88} 
and the properties of the sphere depend crucially on properties of the latter. 

We write  $\A(\CP^n_q)$ for the  $*$-subalgebra of $\A(\rS^{2n+1}_q)$ generated by the elements $p_{ij}:=z_i^*z_j$ for $i,j=0,1,\ldots,n$, which we think of as the coordinate algebra of the quantum projective space $\CP^n_q$. 
It is easy to see that the algebra $\A(\CP^n_q)$ 
is made of the  
invariant elements for the  
action of $\U(1)$ on the algebra $\A(\rS^{2n+1}_q)$ given by
\beq\label{u1act}
(z_0,z_1, \dots, z_n) \mapsto (\lambda z_0, \lambda z_1, \dots, \lambda z_n ), \qquad \lambda \in\U(1).
\eeq
From the relations of $\A(\rS^{2n+1}_q)$ one gets relations for 
$\A(\CP^n_q)$:
\begin{align}\label{rel:qps}
p_{ij}p_{kl}&=q^{\mathrm{sign}(k-i)+\mathrm{sign}(j-l)}\, p_{kl}p_{ij}
 &\hspace{-1.5cm}\mathrm{if}\;i\neq l\;\mathrm{and}\;j\neq k\;, \nn \\
p_{ij}p_{jk}&=q^{\mathrm{sign}(j-i)+\mathrm{sign}(j-k)+1}\, p_{jk}p_{ij}-(1-q^2)
\textstyle{\sum_{l>j}}\, p_{il}p_{lk} &\mathrm{if}\;i\neq k\;, \nn \\
p_{ij}p_{ji}&=
q^{2\mathrm{sign}(j-i)}p_{ji}p_{ij}+
(1-q^2)\left(\textstyle{\sum_{l>i}}\,q^{2\mathrm{sign}(j-i)}p_{jl}p_{lj}
-\textstyle{\sum_{l>j}}\,p_{il}p_{li}\right)
 &\mathrm{if}\;i\neq j\;,
\end{align}
with $\mathrm{sign}(0):=0$.
The elements $p_{ij}$ are  the matrix entries of a projection $P= (p_{ij})$, that is to say it obeys $P^2=P=P^*$, or rather that 
$\sum_{j=0}^n p_{ij} p_{jk} = p_{ik}$ and $p_{ij}^*=p_{ji}$. 
This projection has  $q$-trace equal to one:
\beq\label{q-tr}
\tr_q(P):=\sum_{i=0}^n \, q^{2i} p_{ii}=1.
\eeq
To the best of our knowledge, the algebra $\A(\CP^n_q)$ first appeared in \cite{We00}. 

The $C^*$-algebra $C(\rS^{2n+1}_q)$ of continuous functions on the quantum sphere $\rS^{2n+1}_q$ is the completion of $\A(\rS^{2n+1}_q)$ in the universal $C^*$-norm.
The $C^*$-algebra $C(\CP^n_q)$ of continuous functions on the quantum projective space is the completion of $\A(\CP^n_q)$ in the universal $C^*$-norm. 
By definition, the $*$-algebra inclusion $\A(\CP^n_q)\hookrightarrow \A(\rS^{2n+1}_q)$ extends to an inclusion of $C^*$-algebras $C(\CP^n_q)\hookrightarrow C(\rS^{2n+1}_q)$. 

There is a marked difference between $\rS^{2n+1}_q$ and $\CP^n_q$ which is reflected in their $(K_0,K_1)$-groups:  
for the odd-dimensional spheres $\rS^{2n+1}_q$ these are equal to $(\Z,\Z)$ regardless of the dimension, while for $\CP^n_q$ they are equal to $(\Z^{n+1},0)$. 
A set of generators for the K-theory and K-homology of the sphere algebras $C(\rS^{2n+1}_q)$ can be found in \cite{HL04}.

That $K_0(C(\CP^n_q))\simeq\Z^{n+1}$ can be proved by viewing the $C^*$-algebra $C(\CP^n_q)$ 
as the Cuntz--Krieger algebra of a graph \cite{HS02}. The group $K_0$ is the cokernel of the incidence matrix canonically associated to the graph (while $K_1$ is the kernel of the matrix). The dual result for K-homology is obtained using the same techniques: the group $K^0$ is now the kernel of the transposed matrix \cite{Cu84} and this leads to $K^0(C(\CP^n_q))\simeq\Z^{n+1}$
(and $K^1$ is the cokernel of the transposed matrix).

Generators of the homology group $K^0(C(\CP^n_q))$ were given explicitly in \cite{DL10} as (classes of) even Fredholm modules 
\begin{equation}\label{freds}
\mu_k=(\A(\CP^n_q), \, \HH_{(k)}, \, \pi^{(k)}, \, \gamma_{(k)}, \, F_{(k)}) \,, \quad \text{for} \quad 0\leq k\leq n \,. 
\end{equation}
Generators of the K-theory $K_0(C(\CP^n_q))$  were also given in \cite{DL10} (\textit{cf}. also  \cite{DL13}) as projections whose entries are polynomial functions, that is to say its entries are in the coordinate algebra $\A(\CP^n_q)$ rather than the $C^*$-algebra $C(\CP^n_q)$.

Before we recall these generators explicitly for later use, we need to pause for some notation. 
The $q$-analogue of an integer $n\in \Z$ is  given by
$$
[n]:=\frac{q^n-q^{-n}}{q-q^{-1}} \;;
$$
it is defined for $q\not= 1$ and is equal to $n$ in the limit $q\to 1$. 
For any $n\geq 0$, one defines the factorial of the $q$-number $[n]$ by setting $[0]!:=1$ and then $[n]!:=[n][n-1]\cdots [1]$.  The $q$-multinomial coefficients are in turn defined by
$$
[j_0,\ldots,j_n]!:=\frac{[j_0+\ldots+j_n]!}{[j_0]!\ldots[j_n]!} \;.
$$

For $N\in\Z$, let $\Psi_N:=(\psi^N_{j_0,\ldots,j_n})$ be the vector-valued function on $\rS^{2n+1}_q$ with components
\beq\label{psi}
\psi^N_{j_0,\ldots,j_n}:=
\begin{cases}
\, [j_0,\ldots,j_n]!^{\frac{1}{2}}q^{-\frac{1}{2}\sum_{r<s}j_rj_s} \, (z_0^{j_0})^* \ldots (z_n^{j_n})^* &  \text{for} \quad N \geq 0 \, , \\
& \\
\, [j_0,\ldots,j_n]!^{\frac{1}{2}}q^{\frac{1}{2}\sum_{r<s}j_rj_s+\sum_{r=0}^nrj_r} \, 
z_0^{j_0}\ldots z_n^{j_n} & \text{for} \quad N \leq 0 \, , 
\end{cases}
\eeq
with
$j_0+\ldots+j_n=|N|$. 
Then $\Psi_N^*\Psi_N=1$ and $P_N:=\Psi_N\Psi_N^*$ is a projection in a matrix algebra of a certain size:
\beq\label{size}
P_N\in\textup{M}_{d_N}(\A(\CP^n_q)),\qquad d_N:=\binom{|N|+n}{n} \, ,
\eeq
(this was proven in \cite{DD09}, generalizing the special case where $n=2$ in \cite{DL09}). 
By construction, the entries of the matrix $P_N$ are $\U(1)$-invariant and so they are indeed elements of the algebra $\A(\CP^n_q)$. In particular we see that $P_1=P$ is the `defining' projection of the algebra $\A(\CP^n_q)$ given before with relations in 
\eqref{rel:qps}.  

We let $[P_{N}]$ denote the class in $K_0(C(\CP^n_q))$ of the projection $P_{N}$ and let $[\mu_k]$ denote the classes in $K^0(C(\CP^n_q))$  of the Fredholm modules \eqref{freds}.
The following result was proved in \cite{DL10} (\textit{cf}. Props. 4 and 5 there).

\begin{prop}\label{pair}
For all $N\in\N$ and for all $0\leq k\leq n$ it holds that 
$$
\hp{[\mu_k], [P_{-N}]}:=
\tr_{\HH_k}(\gamma_{(k)} (\pi^{(k)}(\tr\,P_{-N}))
=\tbinom{N}{k} \;,
$$
with $\binom{N}{k}:=0$ when $k>N$. Moreover, the elements $[\mu_0],\ldots,[\mu_n]$  are generators of $K^0(C(\CP^n_q))$, and the elements
$[P_0],\ldots,[P_{-n}]$ are generators of $K_0(C(\CP^n_q))$.
\end{prop}

Indeed, the matrix of couplings $M\in \mathrm{M}_{n+1}(\Z)$ with $M_{ij}:=\hp{[\mu_i],[P_{-j}]}=\tbinom{j}{i}$, for $i,j=0,1,\ldots,n$,
has inverse with integer entries $(M^{-1})_{ij}=(-1)^{i+j}\binom{j}{i}$. Thus the aforementioned elements are a basis of $\Z^{n+1}$ as a $\Z$-module, which is equivalent to saying that they generate $\Z^{n+1}$ as an Abelian group.

\subsection{Line bundles}

It is well known that the algebra inclusion $\A(\CP^n_q) \hookrightarrow \A(\rS^{2n+1}_q)$ 
is a quantum principal bundle with structure group $\U(1)$.
To each projection $P_{N}$ there corresponds a line bundle associated to this principal bundle, as we now describe.   

The column vector $\Psi_N$ has $d_N$ entries, all of which are elements of $\A(\rS^{2n+1}_q)$. 
We consider the collection 
\begin{align}\label{lbcpn}
\LL_N := \left\{ \varphi_N := v \cdot \Psi_N 
=\sum_{j_0+\ldots+j_n=N} v_{j_0,\ldots,j_n} \, \psi^N_{j_0,\ldots,j_n}  \right\} \, ,
\end{align}
where $v=(v_{j_0,\ldots,j_n}) \in (\A(\CP^n_q))^{d_N}$. Each $\LL_N$ is made of elements of $\A(\rS^{2n+1}_q)$ which transform under the $\U(1)$ action in \eqref{u1act}, as $\varphi_N \mapsto \varphi_N \lambda^{-N}$.
In particular $\LL_0=\A(\CP^n_q)$.
By their very definition each $\LL_N$ is an $\A(\CP^n_q)$-bimodule -- the bimodule of equivariant  maps for the irreducible representation of $\U(1)$ with weight $N$. It also holds that
\beq
\LL_N \otimes_{\A(\CP^n_q)} \LL_{N'} \simeq \LL_{N+N'} \,
\eeq
(\textit{cf}. \cite[Lem.~7.5]{DL13} and also \cite[Prop. 3.1]{KLS}) and so, in particular, 
\beq\label{tplb}
(\LL_N)^{\otimes_{\A(\CP^n_q)} M} \simeq \LL_{M N} \,.
\eeq      
An argument as in \cite[Prop.~3.3]{DL10} yields isomorphisms $\LL_{N} \simeq (\A(\CP^n_q))^{d_N}P_{N}$ 
as left $\A(\CP^n_q)$-modules and $\LL_{N}\simeq P_{-N}(\A(\CP^n_q))^{d_N}$ 
as right $\A(\CP^n_q)$-modules. Clearly then, we have to make a choice: we always use the left $\A(\CP^n_q)$-module identification and denote the class of the projection $P_N$  by $[\mathcal{L}_N]$ as an element of the group $K_0(C(\CP^n_q))$.

For each $N\in\Z$ the module $\LL_N$ describes a line bundle, in the sense that its `rank' (as computed by pairing with $[\mu_0]$) is equal to $1$. It is completely characterized by its `first Chern number' (as computed by pairing with the class $[\mu_1]$). Indeed, using an argument similar to that of the proof of Prop.~\ref{pair} one shows the following.

\begin{prop}\label{lb}
For all $N\in\Z$ it holds that 
\begin{align*}
\hp{[\mu_0], [\LL_{N}]} = 1 \, \quad \text{and} \quad 
\hp{[\mu_1], [\LL_{N}]} = -N \, .
\end{align*}
\end{prop}

From the above discussion, the line bundle $\LL_{-1}$ emerges as a central character: from Prop.~\ref{pair} its only non-vanishing charges are $\hp{[\mu_0], [\LL_{-1}]} = 1$ and $\hp{[\mu_1], [\LL_{-1}]} = 1$. The bundle $\LL_{-1}$ is the \emph{tautological line bundle} for the quantum projective space $\CP^n_q$.

Now consider the element in $K_0(C(\CP^n_q))$ given by
\beq\label{ecat}
u := 1 - [\LL_{-1}] \, ,
\eeq
of which we can take powers using the identification \eqref{tplb}. For $j\geq 0$, as elements in K-theory, one has then
\beq
u^j = (1 - [\LL_{-1}])^j \simeq \sum_{N=0}^{j} (-1)^N \tbinom{j}{N} [\LL_{-N}]   \, .
\eeq

\begin{prop}
For $0\leq j \leq n$ and for $0\leq k\leq n$, it holds that
\beq\label{parpow}
\hp{[\mu_k], u^j} = 
\begin{cases}
0 & \text{for} \quad j \not= k  \\
(-1)^j & \text{for} \quad j = k 
\end{cases} \, ,
\eeq
while for all  $0\leq k\leq n$ it holds that
\beq\label{parn1}
\hp{[\mu_k], u^{n+1}} = 0 \,.
\eeq
\begin{proof}
Denoting as before by $[\LL_{-N}]$ the class of the projection $P_{-N}$
and setting $\binom{N}{k}:=0$ when $k>N$, we compute using Prop.~\ref{pair} that
$$
\hp{[\mu_k], u^j} = \sum_{N=0}^{j} (-1)^N \tbinom{j}{N} \hp{[\mu_k],  [\LL_{-N}]} = 
\sum_{N=k}^{j} (-1)^N \tbinom{j}{N} \tbinom{N}{k} \, . 
$$
If $k > j$ this vanishes again due to $\binom{N}{k}:=0$ for $k>N$. On the other hand, if $k \leq j$, it is  
$$
\hp{[\mu_k], u^j} = \frac{j!}{k!} \, \sum_{N=k}^{j} \frac{(-1)^N}{(j-N)!(N-k)!} 
$$
and an act of direct computation yields \eqref{parpow}. Similarly, one computes that 
$$
\hp{[\mu_k], u^{n+1}} = \frac{(n+1)!}{k!} \, \sum_{N=k}^{n+1} \frac{(-1)^N}{(n+1-N)!(N-k)!} = 0 \, ,
$$
thus completing the proof. 
\end{proof}
\end{prop}

The element $u =\chi([\LL_{-1}]) := 1 - [\LL_{-1}]$ shall be named the \emph{Euler class} of the line bundle $\LL_{-1}$, 
in analogy with the classical case (\textit{cf}. \cite[IV.1.13]{Ka78}). 
Since for $0\leq k\leq n$ the elements $[\mu_k]$ are generators of $K^0(C(\CP^n_q))$, the fact that
$\hp{[\mu_k], u^{n+1}} = 0$ for $0\leq k\leq n$ amounts to saying that $u^{n+1}=0$ in $K_0(C(\CP^n_q))$. 
On the other hand, 
since the elements $[\LL_{-k}]$ for $0\leq k\leq n$ are generators of $K_0(C(\CP^n_q))$, the results in \eqref{parpow} say that the elements $[\mu_k]$ and $(-u)^j$ for $0\leq k,j \leq n$ form dual bases. 
These two facts lead to the following analogue of the classical result (\textit{cf}. \cite[Cor.~IV.2.11]{Ka78}). 

\begin{prop}\label{ku}
It holds that
$$
K_0(C(\CP^n_q)) \simeq \Z [\LL_{-1} ] / (1-[\LL_{-1}])^{n+1}  \simeq \Z[u] / u^{n+1} \, 
$$
where $u = \chi([\LL_{-1}]) := 1 - [\LL_{-1}]$ is the Euler class of the line bundle $\LL_{-1}$.
\end{prop}

\section{Quantum lens spaces}\label{qls}

Next we come to describe the algebras of functions on quantum lens spaces and noncommutative `line bundles'  thereon. 
Indeed, we define quantum lens spaces as `direct sums of line bundles' and algebras of their functions in terms of corresponding `modules of sections'. {\em A posteriori} these algebra of functions are seen as subalgebras of functions on 
odd-dimensional quantum spheres which are invariant for the action of a cyclic group. 
In such a manner there are natural principal and associated fibrations. 
This gives rise to a natural family of representatives of classes in the K-theory of quantum lens spaces.

\subsection{Functions on quantum lens spaces} Fix an integer $r \geq 2$ and define
\beq\label{dals}
\A(\rL^{(n,r)}_q) :=\bigoplus_{N\in\Z} \, \LL_{r N} \, .
\eeq
Then proving the following is straightforward.  

\begin{prop}\label{zract}
The vector space $\A(\rL^{(n,r)}_q)$ is a $*$-algebra made of all elements of $\A(\rS^{2n+1}_q)$ which are invariant under the action $\alpha_r :\Z_r\to \textup{Aut}(\A(\rS^{2n+1}_q))$ of the cyclic group $\Z_r$ generated by the map
\beq \label{gpaction}
(z_0, z_1, \dots, z_n) \mapsto ( e^{ 2 \pi \ii / r} z_0, e^{ 2 \pi \ii / r} z_1, \dots, e^{ 2 \pi \ii / r} z_n ) \,.
\eeq 
\end{prop}
\proof 
Let $\zeta = e^{ 2 \pi \ii / r}$, which then satisfies $\zeta^r=1$, be the generator of $\Z_r$. 
A standard argument shows that equipping the algebra $\A(\rS^{2n+1}_q)$ with this action $\alpha_r$ 
is equivalent to giving it the grading by the Pontryagin dual group $\widehat \Z_r \simeq \Z_r$ defined by
$$
\A(\rS^{2n+1}_q)=\bigoplus_{k\in\Z_r} \A_k,\qquad \A_k:=\{a\in \A(\rS^{2n+1}_q)~|~\alpha_r(a) = \zeta^{-k} \, a\}.
$$
The algebra $\A_0$ of invariant elements is by definition the algebra $\A (\rL^{(n,r)}_q)$ of functions on the lens space; 
this  also shows that there are no other invariant elements. 
\endproof

We think of the algebra $\A(\rL^{(n,r)}_q)$ as the coordinate algebra of an underlying quantum space $\rL^{(n,r)}_q$, 
which is named the \textit{quantum lens space} of dimension $2n+1$ (and index $r$); it is a deformation of the classical lens space 
$L^{(n,r)}= \rS^{2n+1} / \Z_r$ of the same dimension. 
The $C^*$-algebra $C(\rL^{(n,r)}_q)$ of continuous functions on the quantum lens space, the universal $C^*$-completion 
of $\A(\rL^{(n,r)}_q)$, is part of the general family of lens spaces defined in \cite{HS03}.
 
Of course, the value $r=1$ is also possible but this does not yield anything new. Indeed, in that case one 
has $L^{(n,1)}_q = \rS^{2n+1}_q$ and the above expression \eqref{dals} is nothing other than the well known
vector space decomposition 
\beq\label{Lnbdl}
 \A (\rS^{2n+1}_q) = \bigoplus_{N\in\Z} \, \LL_{N} \, .
\eeq
Clearly $\A(\CP^n_q)$ is a subalgebra of $\A(\rL^{(n,r)}_q)$. In parallel with the $\U(1)$ quantum principal bundle $\A(\CP^n_q) \hookrightarrow \A(\rS^{2n+1}_q)$ there is indeed more structure. 

\begin{prop}\label{pr-u1t}
The algebra inclusion $\A(\CP^n_q) \hookrightarrow \A(\rL^{(n,r)}_q)$ is a quantum principal bundle with structure group  
$\widetilde{\U}(1) :=\U(1)/\Z_r \simeq \U(1)$.
In particular, one finds that
$$
\A(\CP^n_q) = \A(\rL^{(n,r)}_q)^{\widetilde{\U}(1)},
$$
in analogy with the identification $\A(\CP^n_q) = \A(\rS^{2n+1}_q)^{\U(1)}$ as defined before.
\end{prop}
\proof This is deferred to  App.~\ref{ncbundles}. 
\endproof

The $\U(1)$ and $\widetilde\U(1)$ principal bundles over the quantum projective space
$\CP^n_q$ as in Prop.~\ref{pr-u1t} are related by a $\Z_r$ principal bundle structure over the quantum lens space $\rL^{(n,r)}_q$. 
Indeed, it is also not difficult to verify the following.
\begin{prop}\label{pr-Zr}
The algebra inclusion $\A(\rL^{(n,r)}_q) \hookrightarrow \A(\rS^{2n+1}_q)$ is a quantum principal bundle with structure group $\Z_r$.  
\end{prop}
\proof This can also be found in  App.~\ref{ncbundles}.
\endproof

\begin{rem}\textup{
Note that, with $\widetilde\Z$ denoting the Pontryagin dual of the quotient group $\widetilde\U(1)$ (so that there is an injection $\Z\hookrightarrow \widetilde\Z$ given by multiplication by $r$), the decomposition \eqref{dals} may be written, in parallel with the decomposition \eqref{Lnbdl}, as
\beq\label{dals2}
\A(\rL^{(n,r)}_q)=\bigoplus_{N\in \widetilde\Z}\LL_N .
\eeq
}
\end{rem}

\subsection{Pulling back line bundles}

Having the principal bundle $j: \A(\CP^n_q) \hookrightarrow \A(\rL^{(n,r)}_q)$, we proceed now to `pull-back' 
the associated line bundles from $\CP^n_q$ to $\rL^{(n,r)}_q$. We are led to the following natural definition.

\begin{defi}\label{pblb}
For each $\A(\CP^n_q)$-bimodule $\LL_N$ as in \eqref{lbcpn}  (a line bundle over $\CP^n_q$),  its `pull-back' to $\rL^{(n,r)}_q$ is the 
$\A(\rL^{(n,r)}_q)$-bimodule
\beq
j_*(\LL_N) := \left\{ \widetilde{\varphi}_N = v \cdot \Psi_N =\sum_{j_0+\ldots+j_n=N} v_{j_0,\ldots,j_n} \, \psi^N_{j_0,\ldots,j_n} \right\} \, ,
\eeq
for $v=(v_{j_0,\ldots,j_n}) \in (\A(\rL^{(n,r)}_q))^{d_N}$. We shall often use the shorthand $j_*(\LL_N) :=   \widetilde{\LL}_N$.
\end{defi} 

By embedding the cyclic group $\Z_r$ into $\U(1)$ via the $r$-th roots of unity, each $\widetilde{\LL}_N$  is made of elements of $\A(\rS^{2n+1}_q)$ which transform as $\widetilde{\varphi}_N \mapsto \widetilde{\varphi}_N \, e^{ - 2 \pi \ii \, N / r} $ under the $\U(1)$-action 
of Prop.~\ref{zract}. By its very definition, $\widetilde{\LL}_N$ is an $\A(\rL^{(n,r)}_q)$-bimodule. Once again, arguments like those of \cite[Prop.~3.3]{DL10} for the $\widetilde{\LL}_N$ yield the following.

\begin{prop}
There are left $\A(\rL^{(n,r)}_q)$-module isomorphisms
$$
\widetilde{\LL}_{N} \simeq (\A(\rL^{(n,r)}_q))^{d_N} P_{N}   
$$
and right $\A(\rL^{(n,r)}_q)$-module isomorphisms
$$
\widetilde{\LL}_{N}\simeq P_{-N}(\A(\rL^{(n,r)}_q))^{d_N} \, .
$$ 
\end{prop}

We stress that the projections $P_N$ here are those constructed before, around \eqref{psi} and \eqref{size},
taken now as elements of the group $K_0(C(\rL^{(n,r)}_q))$.
Just as for the modules $\LL_{N}$, we need to make a choice of representatives: we use the left $\A(\rL^{(n,r)}_q)$-module identification and denote by $[\widetilde{\LL}_{N}]$ the class of the projection $P_{N}$ as an element in $K_0(C(\rL^{(n,r)}_q))$.
Thus, the pull-back of line bundles induces a map
\beq\label{Kpblb}
j_* : K_0(C(\CP^n_q)) \to K_0(C(\rL^{(n,r)}_q)) \, .
\eeq
`Geometrically', the pull-back of line bundles from $\CP^n_q$ to $\rL^{(n,r)}_q$ could be depicted as   
\beq\label{dia-pb}
\xymatrix{ 
\widetilde{\LL}_N \ar@{.>}[d] & \ar[l]_{j_{*}} \LL_N \ar@{.>}[d] & \\  
 \A(\rL^{(n,r)_q}) \,\, & \,\, \ar[l]^j \A(\CP^n_q) \, .  &
}
\eeq

The marked difference between the module $\LL_N$ versus its pull-back $\widetilde{\LL}_{N}$ is that, while each $\LL_N$ 
is not free when $N\ne 0$ (as a consequence of Prop.~\ref{lb}), this need not be the case for  
$\widetilde{\LL}_{N}$, that is the projection $P_N$ could be trivial (i.e. equivalent to $1$) in $K_0(C(\rL^{(n,r)}_q))$.

Indeed, the pull-back $\widetilde{\LL}_{-r}$ of the line bundle $\LL_{-r}$  
from the projective space $\CP^n_q$ to the lens space $\rL^{(n,r)}_q$ is free: recall that the corresponding projection is $P_{-r}:=\Psi_{-r}\Psi_{-r}^*$ and here the vector-valued function $\Psi_{-r}$ has entries in the algebra $\A(\rL^{(n,r)}_q)$ itself. Thus the condition $\Psi_{-r}^*\Psi_{-r}=1$ implies that the projector $P_{-r}$ is equivalent to $1$, that is to say,
the class of the module $\widetilde{\LL}_{-r}$ is trivial in the group $K_0(C(\rL^{(n,r)}_q))$. 
It follows that $(\widetilde{\LL}_{-N})^{\otimes r}\simeq\widetilde{\LL}_{-rN}$ also has trivial class for any $N\in\Z$, the tensor product being taken over $\A(\rL^{(n,r)}_q)$.  Such pulled-back line bundles $\widetilde{\LL}_{-N}$ thus define \emph{torsion classes} and, as we shall see later on, they generate the group $K_0(C(\rL^{(n,r)}_q))$. 

\begin{rem}\textup{
A moment's thought shows that for each $\A(\CP^n_q)$-module $\LL_N$ its pull-back  $\widetilde{\LL}_N$ 
is none other than the $\A(\rL^{(n,r)}_q)$-module
$$
\widetilde{\LL}_N = \A(\rL^{(n,r)}_q) \otimes_{\A(\CP^n_q)} \LL_N \, .
$$
From this it follows at once that $\widetilde{\LL}_{-r} = \A(\rL^{(n,r)}_q) \otimes_{\A(\CP^n_q)} \LL_{-r} \simeq \A(\rL^{(n,r)}_q) = \widetilde{\LL}_{0} 
$, thus showing that the module $\widetilde{\LL}_{-r}$ is free. While we could have taken this as defining the pull-back map, we rather prefer the one in \eqref{pblb} due to the central role to be played by the partial isometries $\Psi_N$'s later on in the paper.
}
\end{rem}

At this point it is pertinent to introduce a second crucial ingredient to the discussion, in the form of a natural map
\beq\label{Gsqls}
\alpha : K_0(C(\CP^n_q)) \to K_0(C(\CP^n_q)),
\eeq
where $\alpha$ is multiplication by the Euler class $\chi(\LL_{-r}) := 1 - [\LL_{-r}]$ of the line bundle $\LL_{-r}$. The central idea of our paper is to combine this map with the pull-back map \eqref{Kpblb} into a sequence for our quantum lens spaces that parallels the classical Gysin sequence \eqref{cgs}:
$$
\xymatrix{ 
0\to K_1(C(\rL^{(n,r)}_q)) \ar[r]^-{\textup{Ind}} & K_0(C(\CP^n_q)) \ar[r]^-{\alpha} 
& K_0(C(\CP^n_q)) \ar[r]^-{j_*} & K_0(C(\rL^{(n,r)}_q))  \ar[r] & 0\, ,
}
$$
for a suitable map $\textup{Ind}$. 
Such is the topic of the next section.

\section{The Gysin sequence for quantum lens spaces}

In this section we arrive at a Gysin sequence for quantum lens spaces by invoking some general properties of K-theory associated to circle actions on $C^*$-algebras. For the general background of unbounded operators on Hilbert modules we refer to \cite{Lance}. 
Given a (countably generated) right Hilbert $F$-module $X\leftrightharpoons F$ over a $C^*$-algebra $F$, 
the right $F$-Hermitian structure on $X$ is denoted $\langle\cdot |\cdot\rangle_F$. 

\subsection{Construction of the sequence} \label{cose}
In order to lighten our notation and in a way which is consistent with the notation used in \cite{CPR, CNNR08}, we write
$$
A:=C(\rL^{(n,r)}_q),\qquad F:=C(\CP^n_q) .
$$ 

By the universal property of these $C^*$-algebras, the action of $\widetilde\U(1)$ on $\A(\rL^{(n,r)}_q)$ extends uniquely to a strongly continuous circle action $\sigma:\widetilde\U(1)\to \textup{Aut}(A)$ on $A$. From Prop.~\ref{pr-u1t}
the $C^*$-algebra $F$ sits inside $A$ as the fixed point subalgebra, namely 
$$
F=\{a\in A : \sigma_t(a)=a ~ \text{for all}~t\in \widetilde\U(1)\}.
$$
Since $\widetilde\U(1)$ is compact, there is a faithful conditional expectation
$$
\tau:A\to F,\qquad \tau(a):=\int_0^{2\pi} \sigma_t(a)\textup{d} t \,,
$$
from which one obtains an $F$-valued inner product on $A$ by defining
$$
\langle\cdot,\cdot\rangle_F: A\times A\to F,\qquad \langle a,b\rangle_F:=\tau(a^*b).
$$
The properties of the conditional expectation imply that this equips $A$ with the structure of a right pre-Hilbert $F$-module. Let $X\leftrightharpoons F$ be the right Hilbert module resulting from completion of $A$ in the corresponding norm $\|a\|_X^2:=\|\langle a,a\rangle_F\|$. For each $k\in \Z$ we denote the corresponding eigenspace of the action $\sigma$ on $A$ by
$$
A_k:=\{ a\in A : \sigma_t(a)=e^{i k t } a  ~\text{for all} ~t\in \R\}.
$$
In particular $A_0=F$ and in fact Prop.~\ref{pr-u1t} implies that the action  $\sigma:\widetilde\U(1)\to \textup{Aut}(A)$ has {\em full spectral subspaces} in the sense that  $\overline{A_k^*A_k}=F$ for all $k\in\Z$ (cf. \cite[Defn.~2.2]{CNNR08}). This is interpreted as a noncommutative analogue of having a free circle action.

We recall that if $\mathfrak{D}:\Dom (\mathfrak{D})\rightarrow X$ is an unbounded linear operator with dense domain $\Dom(\mathfrak{D})\subseteq X$, then $\mathfrak{D}$ is said to be {\em closed} whenever its graph
\begin{equation}
\mathfrak{G}(\mathfrak{D}):=\left\{\begin{pmatrix}x, \mathfrak{D}x\end{pmatrix}~|~x\in \Dom (\mathfrak{D})\right\}\subseteq X\oplus X
\end{equation}
is a closed subspace of $X\oplus X$. The operator $\mathfrak{D}:\Dom(\mathfrak{D})\to X$ is said to be {\em symmetric} if $\Dom(\mathfrak{D})\subseteq \Dom(\mathfrak{D}^*)$ and $\mathfrak{D}=\mathfrak{D}^*$ on $\Dom(\mathfrak{D})$, and it is {\em self-adjoint} if it is symmetric and $\Dom(\mathfrak{D}^*)=\Dom(\mathfrak{D})$. A closed, self-adjoint linear operator $\mathfrak{D}:\Dom (\mathfrak{D})\rightarrow X$ is said to be {\em regular} if and only if the operators $\mathfrak{D}\pm i:\Dom (\mathfrak{D})\rightarrow X$ have dense range, which in turn happens if and only if these operators are bijective.

There is a certain natural unbounded operator with these properties on the Hilbert module $X\leftrightharpoons F$ which is of particular interest to us. This operator is nothing other than the infinitesimal generator of the circle action  
$\sigma:\widetilde\U(1)\to \textup{Aut}(A)$ given explicitly by
\beq\label{suad}
\mathfrak{D}:X_\mathfrak{D}\to X,\qquad \mathfrak{D} \left(\sum_{k\in\Z} x_k\right):=\sum_{k\in\Z} k \, x_k
\eeq
on the dense domain $X_\mathfrak{D}\subset X$,
\beq\label{suax}
X_\mathfrak{D}:=\left\{ x=\sum_{k\in\Z} x_k \in X~|~ x_k\in A_k,~\left\| \sum_{k_\in\Z}k^2\langle x_k,x_k\rangle\right\|<\infty\right\} . 
\eeq
This defines a self-adjoint and regular operator on $X$ (\textit{cf}. \cite[Prop.~4.6]{PR} and \cite[Prop.~2.7]{CNNR08}) and
it follows from  \cite[Prop.~2.9]{CNNR08} that the pair $(X,\mathfrak{D})$ in turn yields a class in the odd unbounded Kasparov bivariant K-theory $KK_1(A,F)$. 
 
Next, we invoke the internal Kasparov product in bivariant K-theory, which in our case of interest takes the form of a map
$$
-\,\widehat\otimes_A-:KK_*(\C,A)\times KK_1(A,F)\to KK_{*+1}(\C,F).
$$
In particular, the Kasparov product of the class $[(X,\mathfrak{D})]$ of the pair $(X,\mathfrak{D})$ in $KK_1(A,F)$ with the K-theory $K_*(A)=KK_*(\C,A)$ immediately equips us with a pair of maps
$$
\textup{Ind}_\mathfrak{D}:K_*(A)\to K_{*+1}(F) , \qquad \textup{Ind}_\mathfrak{D}(-):=-\,\widehat\otimes_A [(X,\mathfrak{D})].
$$

For the case of interest in the present paper one has $K_1(F)=0$, thus one of these maps is just the zero map, $\textup{Ind}_\mathfrak{D}: K_0(A)\to 0$. On the other hand, 
by its definition the operator $\mathfrak{D}$ has a spectral gap around zero, whence the map $\textup{Ind}_\mathfrak{D}:K_1(A)\to K_0(F)$ 
is  given explicitly as an index
$$
\textup{Ind}_\mathfrak{D}([u]):=[\textup{Ker}\,PuP]-[\textup{Coker}\,PuP]\in K_0(F),
$$
where $P$ denotes the spectral projection for the self-adjoint operator $\mathfrak{D}$ and associated to the non-negative real axis \cite[App.~A]{PR}.

In our new notation, the multiplication in \eqref{Gsqls} by the Euler class $\chi(\LL_{-r}) := 1 - [\LL_{-r}]$ of the line bundle $\LL_{-r}$ yields a map
$$
\alpha : K_0(F) \to K_0(F),
$$
whereas the map in \eqref{Kpblb} induced by the inclusion map $j:F\hookrightarrow A$ gives
$$
j_* : K_0(F) \to K_0(A).
$$
Assembling all of these together yields a sequence
\begin{equation}\label{qgysin}
\xymatrix{ 
0\to K_1(A) \ar[r]^-{\textup{Ind}_\mathfrak{D}} & K_0(F) \ar[r]^-{\alpha} & K_0(F) \ar[r]^-{j_*} & K_0(A) \to 0\, ,
}
\end{equation}
where we already use the fact that in the last term of the sequence we have $K_1(F)=0$. We claim that it is an {\em exact} sequence: 
this is proved over the next couple of sections. The sequence \eqref{qgysin} will be called the {\em Gysin sequence} for the quantum lens space $\rL^{(n,r)}_q$. 

\begin{rem}\textup{
Let us stress that it is not invidious to put zeros at the beginning and end of the sequence \eqref{qgysin}. At this point we are saying nothing about exactness of the sequence, so we are not yet saying that the maps $\textup{Ind}_{\mathfrak{D}}$ and $j_*$ are respectively injective and surjective: for the time being this is merely a claim.
}
\end{rem}

\subsection{K-theory of the mapping cone} 
Our strategy to prove exactness of the Gysin sequence will be to relate it to a six-term exact sequence in K-theory coming from the mapping cone of the pair $(F,A)$. We start then by recalling how to lift the index theory described above to the KK-theory of the mapping cone construction.
 
 Recall that the {\em mapping cone} of the pair $(F,A)$ is the $C^*$-algebra
$$
M(F,A):=\left\{ f\in C([0,1],A)~|~ f(0)=0,~ f(1)\in F \right\}.
$$
The group $K_0(M(F,A))$ has a particularly elegant description in terms of partial isometries. Indeed, let us write $V_m(F,A)$ for the set of partial isometries $v\in \textup{M}_m(A)$ such that the associated projections $v^*v$ and $vv^*$ belong to $\textup{M}_m(F)$. Using the inclusion $V_m(F,A)\hookrightarrow V_{m+1}(F,A)$ given by setting $v\mapsto v\oplus 0$, one defines
$$
V(F,A):=\bigcup_m V_m(F,A)
$$
and then generates an equivalence relation $\sim$ on $V(F,A)$ by declaring that:
\begin{enumerate}
\item $v\sim v\oplus p$ for all $v\in V(F,A)$ and $p\in \textup{M}_m(F)$; \\
\item if $v(t)$, $t\in [0,1]$, is a continuous path in $V(F,A)$ then $v(0)\sim v(1)$.
\end{enumerate}
Following \cite[Lem.~2.5]{IP}, there is a well defined bijection between $V(F,A)/\!\!\sim$ and  the K-theory $K_0(M(F,A))$. It is also shown there that if $v$ and $w$ are partial isometries with the same image in 
$K_0(M(F,A))$ one can arrange them to have the same initial projection, i.e. $v^*v=w^*w$, without changing their class in 
$V(F,A)/\!\!\sim$. Having done so, an addition is defined \cite[Lem.~3.3]{CPR} in $V(F,A)/\!\!\sim$ by  
$[v\oplus w^*]=[v]+[w^*]=[v]-[w]=[vw^*]$ so that $V(F,A)/\!\!\sim$ and $K_0(M(F,A))$ 
are isomorphic as Abelian groups.
 
The class of $(X,\mathfrak{D})$ in $KK_1(A,F)$, as defined previously in \S\ref{cose}, has a canonical lift to the group $KK_0(M(F,A),F)$. 
Let $P$ be as before, the spectral projection for $\mathfrak{D}$ corresponding to the non-negative real axis. Following \cite[\S 4]{CPR}, 
we write 
$$
T_{\pm}:=\pm\partial_t\otimes 1+1\otimes \mathfrak{D}
$$
for the unbounded operators with domains
\begin{align*}
\mathfrak{Dom}(T_\pm):&=\left\{ f\in C^\infty_c([0,\infty))\otimes X_\mathfrak{D}~|~ f=\sum_{i=1}^n f_i\otimes x_i,~ x_i\in X_\mathfrak{D}, \right\} \\ &\textup{and} \qquad\quad P(f(0))=0 ~\text{(+ case)}, \quad (1-P)(f(0))=0 ~\text{(-- case)} ,
\end{align*}
where smoothness at the boundary of $[0,\infty)$ is defined by taking one-sided limits. 
With $Y:=L^2([0,\infty))\otimes X$, one finds that
$$
\widehat{\mathfrak{D}}:\mathfrak{Dom}(T_+)\oplus\Dom(T_-)\to Y \oplus Y ,
\qquad \widehat{\mathfrak{D}}:=\begin{pmatrix} 0 & T_- \\ T_+ & 0 \end{pmatrix} ,
$$
is a densely defined unbounded symmetric linear operator. By modifying the domains slightly, one obtains a $\Z_2$-graded Hilbert $M(F,A)$-$F$-bimodule $\widehat X$ endowed with an odd unbounded linear operator $\widehat{\mathfrak{D}}:\Dom(\widehat{\mathfrak{D}})\to\widehat X$ which in addition is self-adjoint and regular \cite[Prop.~4.13]{CPR}. It then follows from  \cite[Prop.~4.14]{CPR} that the pair $(\widehat{X},\widehat{\mathfrak{D}})$ determines a class in the bivariant K-theory $KK_0(M(F,A),F)$.

The internal Kasparov product of $K_0(M(F,A))$ with the class of $(\widehat{X},\widehat{\mathfrak{D}})$ yields a map
\begin{equation}\label{ind lift}
\textup{Ind}_{\widehat{\mathfrak{D}}} : K_0(M(F,A))\to K_0(F).
\end{equation}
Following \cite[Thm.~5.1]{CPR} and writing $\mathcal{E}:=X^m$, the internal Kasparov product of the K-theory $K_0(M(F,A))$ with the class of $(\widehat{X},\widehat{\mathfrak{D}})$ in the bivariant K-theory $KK_0(M(F,A),F)$ is represented by the index
$$
\textup{Ind}_{\widehat{\mathfrak{D}}}([v]):=\textup{Ker}(PvP)|_{v^*vP\mathcal{E}} - \textup{Ker}(Pv^*P)|_{vv^*P\mathcal{E}},
$$
the result being an element of $KK_0(\C,F)=K_0(F)$. 
Here $v\in \textup{M}_m(A)$ is a partial isometry representing a class in $K_0(M(F,A))$ and considered as a map $v: v^*vP\mathcal{E}\to vv^*P\mathcal{E}$.

\subsection{Exactness of the Gysin sequence}  With $S(A):=C_0((0,1))\otimes A$ the {\em suspension} of the $C^*$-algebra $A$, it is clear that 
\begin{equation}\label{susex}
0\to S(A) \xrightarrow{i} M(F,A) \xrightarrow{\textup{ev}} F \to 0,
\end{equation}
where $i(f\otimes a)(t):=f(t)a$ and $\textup{ev}(f):=f(1)$, is an exact sequence of $C^*$-algebras. The six term exact sequence in K-theory corresponding to \eqref{susex} has the form
\begin{equation}\label{six2}
\xymatrix{ 
K_0(S(A))\ar[r] & K_0(M(F,A))\ar[r] & K_0(F)\ar[d] \\  K_{1}(F)\ar[u] & K_{1}(M(F,A))\ar[l] & K_{1}(S(A)). \ar[l]
}
\end{equation}
Using the vanishing of $K_1(F)$, this sequence degenerates to 
\begin{equation}\label{six3}
\xymatrix{
0 \to K_0(S(A)) \ar[r]^-{i_*} & K_0(M(F,A)) \ar[r]^-{\textup{ev}_*} & K_0(F) \ar[r]^-{\partial} & K_1(S(A)) \to K_1(M(F,A))\to 0\, .
}
\end{equation}

As before, the inclusion $i:S(A)\to M(F,A)$ induces the map 
$$
i_*:K_0(S(A))\to K_{0}(M(F,A)).
$$
The map $\textup{ev}_*$ in \eqref{six3} can be given by 
$$
\textup{ev}_* : K_0(M(F,A))\to K_0(F),\qquad \textup{ev}_*([v]):=[v^*v]-[vv^*],
$$
for $v\in \textup{M}_m(A)$ a partial isometry representing a class in $K_0(M(F,A))$ (\emph{cf}. \cite[Lem.~2.3]{IP}). 

The boundary map $\partial$ is defined as in \cite[p.113]{HR}: for $[p]-[q]\in K_0(F)$ one chooses representatives $p,q$ over $F$ and, from these, self-adjoint lifts $x,y$ over $M(F,A)$. Then the exponentials $e^{2\pi i x}$ and $e^{2\pi i y}$ are unitaries over $C(S^1)\otimes A$ which are equal to the identity modulo $C_0((0,1))\otimes A$, so one defines
\begin{equation}\label{exp map}
\partial([p]-[q]):=[e^{2\pi i x}] -[e^{2\pi i y}]\in K_1(S(A)).
\end{equation}

Now recall that $A$ is a Cuntz-Krieger algebra associated to a graph which is connected, row-finite and has neither sources nor sinks  \cite{HS03}. It follows  \cite[Lem.~6.7]{CPR} that $K_1(M(F,A))=0$ and that the index map \eqref{ind lift} is an isomorphism 
\cite[Prop.~6.8]{CPR}. Thus 
\beq\label{komfa}
K_0(M(F,A)) \simeq K_0(F) \simeq \Z^{n+1} ,
\eeq
where the second isomorphism is the result of Prop.~\ref{pair}, since $F=C(\CP^n_q)$.

As a consequence, there is a very easy description of the partial isometries which generate $K_0(M(F,A))$. Recall from the discussion at the end of \S\ref{qls}   that, upon pulling-back to $A$, one finds that for all $N\in \Z$ the projections $P_{rN}$ become equivalent to the identity. This is equivalent to saying that for any $M\in \Z$ the projections $P_{rN}$ and $P_{r(N+M)}$ are equivalent for all $N\in\Z$. Indeed, one can explicitly exhibit partial isometries relating these projectors. Taking the particular case $M=1$, these partial isometries are the
elements $v_N\in \textup{M}_{(d_{r(N+1)},\, d_{rN})}(A)$, with the integers $d_{(\cdot)}$ as in \eqref{size}, given by 
\beq\label{partials}
v_N = \Psi_{r(N+1)} \, \Psi_{rN}^\dag, \qquad N=0,-1,\ldots, -n ;
\eeq
clearly $v_N^*v_N=P_{rN}$ and $v_Nv_N^*=P_{r(N+1)}$ for $N=0,-1,\ldots, -n$. 
With our conventions,  the entries of $v_N$ are elements of $A$ homogeneous of degree $-r$ for the action 
of $\widetilde \U(1)$.

\begin{prop}
The partial isometries \eqref{partials} form a basis of $K_0(M(F,A))$.
\end{prop}

\proof 
From \eqref{komfa} we just need $n+1$ independent generators. Now, since the map \eqref{ind lift} is an isomorphism, the partial isometries $v_N$ are independent (and thus a basis for $K_0(M(F,A)$) if and only if the classes $\textup{Ind}_{\widehat{\mathfrak{D}}}([v_N])$ are so. Since $Pv_NP$ is essentially a `left degree shift' operator on the elements of non-negative homogeneous degree in $PX^{d_{rN}}$ it has no cokernel. Its kernel thus determines the index:
$$
\textup{Ind}_{\widehat{\mathfrak{D}}}([v_N])=\left[P_{rN}X_0^{d_{rN}}\right]=[P_{rN}].
$$
Now, it follows from Prop.~\ref{pair} that the matrix of pairings $\{ \hp{[\mu_k], [P_{-rN}]} = \binom{rN}{k} \}$ is invertible, thus proving that the elements $P_{rN}$
for $N=0,-1,\ldots, -n$ are independent. We note that these projections do not form a basis for $K_0(F)$: the matrix of pairings (while invertible over $\mathbb{Q}$) is not invertible over $\Z$, that is it does not belong to $\textup{GL}(n+1,\Z)$. 
\endproof

Finally we introduce a pair of maps 
$$
B_F:K_0(F)\to K_0(F),\qquad B_A:K_1(S(A))\to K_0(A).
$$
The former is defined simply by the multiplication 
$$
B_F([p]-[q]):=-[P_{-r}]([p]-[q])
$$
in $K_0(F)$. The latter map $B_A$ is the inverse of the Bott isomorphism 
$$
\textup{Bott}:K_0(A)\to K_1(S(A))
$$ 
given by $\textup{Bott}([p]):=[e^{-2\pi i t}\otimes p + 1\otimes (1-p)]$, with $t\in (0,1)$, where $S(A)=C_0((0,1))\otimes A$. 

We are ready to state and prove\footnote{We thank Adam Rennie for explaining to us the explicit forms of the various maps in the exact sequence \eqref{six3}, a conversation from which the proof of our theorem followed very naturally.} a central theorem; it 
directly implies exactness of the Gysin sequence \eqref{qgysin}.

\begin{thm}
There is a diagram
\begin{equation}
\xymatrix{
0 \ar[r] & K_1(A) \ar[r]^-{i_*} \ar[d]^{\textup{id}} & K_0(M(F,A)) \ar[r]^-{\textup{ev}_*} \ar[d]^-{\textup{Ind}_{\widehat{\mathfrak{D}}}} & K_0(F) \ar[r]^-{\partial}  \ar[d]^-{B_F} & K_1(S(A))  \ar[r] \ar[d]^-{B_A} & 0 \\
0\ar[r]  & K_1(A) \ar[r]^-{\textup{Ind}_\mathfrak{D}} & K_0(F) \ar[r]^-{\alpha} & K_0(F) \ar[r]^-{j_*} & K_0(A)  \ar[r] & 0
}
\end{equation}
in which every square commutes and each vertical arrow is an isomorphism of groups. 
\end{thm}

\proof 
Upon using the isomorphism $K_1(A)\simeq K_0(S(A))$, that the first square commutes is precisely \cite[Thm.~5.1]{CPR}.  For the second square we explicitly compute that for each 
$N=0,-1,\ldots, -n$, one has
\begin{align*}
\alpha\left(\textup{Ind}_{\widehat{\mathfrak{D}}}([v_N])\right)=(1-[P_{-r}])[P_{rN}]
=-[P_{-r}]([P_{rN}]-[P_{r(N+1)}])=B_F(\textup{ev}_*([v_N]).
\end{align*}
For the third square, we argue as in \cite[Lem.~3.1]{CPR}.  Recall that in defining the map \eqref{exp map} we chose self-adjoint lifts $x,y$ over $M(F,A)$. We choose here in particular the lifts $x:=t\otimes j(p)$ and $y:=t\otimes j(q)$. These are both self-adjoint and vanish at $t=0$; at $t=1$ they are matrices over $F$. It follows that
$$
[e^{2\pi i x}] -[e^{2\pi i y}]=[e^{2\pi i(t\otimes p)}] -[e^{2\pi i (t\otimes q)}]=-\textup{Bott}([p]-[q])\in K_1(S(A)).
$$
Thus it follows that, modulo the isomorphism $\textup{Bott}:K_0(A)\to K_1(S(A))$, we have 
\begin{equation}\label{exp-j}
\partial([p]-[q])=-([j(p)]-[j(q)]),
\end{equation}
i.e. that $\partial$ is induced up to Bott peridocity by minus the algebra inclusion $j:F\to A$. Now using the fact that the image of the class of $P_{-r}$ in $K_0(A)$ along $j:F\to A$ is trivial, the above \eqref{exp-j} may in fact be written
$$
\partial([p]-[q])=-[j(P_{-r})]([j(p)]-[j(q)]) ,
$$
up to Bott periodicity, from which the result follows.\endproof

\section{The K-theory of quantum lens spaces} \label{ktqls}

We put to work the Gysin sequence \eqref{qgysin} by using it  to compute the K-theory of our quantum lens spaces. 
We shall obtain explicit generators as classes of `line bundles', generically torsion ones.  
This is illustrated by working out some explicit examples.

Since the map $j_*$ in \eqref{qgysin} is surjective, the group $K_0(C(\rL^{(n,r)}_q))$ 
can be obtained by `pulling back' classes from $K_0(C(\CP^n_q))$. Now, 
as shown in Prop.~\ref{ku}, 
$$
K_0(C(\CP^n_q)) \cong \Z[u] / u^{n+1} 
$$ 
with $u := 1 - [\LL_{-1}]$ the Euler class of the line bundle $\LL_{-1}$.  
Moreover,  the Euler class $\chi(\LL_{-r})$ of the line bundle $\LL_{-r}$ is just 
$$
\chi(\LL_{-r}) = 1 - [\LL_{-r}] = 1 - [\LL_{-1}]^r = 1 - (1-u)^r .
$$ 
As a consequence, the map $\alpha$ in \eqref{qgysin} can be given 
as an $(n+1)\times (n+1)$ matrix $A$ with respect to the $\Z$-module basis 
$\{1, u, \dots, u^n \}$ of $K^0(C(\CP^n_q)) \simeq \Z^{n+1}$. This leads to
\beq\label{qk1ls}
K_1(C(\rL^{(n,r)}_q)) \simeq \ker(\alpha) = \ker(A), \qquad K_0(C(\rL^{(n,r)}_q)) \simeq \coker(\alpha) = \coker(A)  \, , 
\eeq
as $\Z$-module identifications via the surjective `pull-back' map $j_*$. 

Simple algebra allows one to compute explicitly the matrix $A$ 
of the map $\alpha$ with respect to the $\Z$-module basis 
$\{1, u, \dots, u^n \}$. Using the condition $u^{n+1}=0$ one has
$$
\chi(\LL_{- r}) = 1-(1-u)^r = \sum_{j=1}^{\min{(r,n)}} (-1)^{j+1} \tbinom{r}{j} u^j \, .
$$
Thus $A$ is an $(n+1) \times (n+1)$ strictly lower triangular matrix
with entries on the $j$-th sub-diagonal equal to $(-1)^{j+1} \tbinom{r}{j}$ for $j \leq \min({r,n})$ and zero otherwise:
\beq\label{matalp}
A= \begin{pmatrix}
0 & 0&0 & \cdots \cdots & & & 0\\
r & 0 & 0 & \cdots \cdots & & & 0 \\
-\tbinom{r}{2} & r & 0 & \cdots \cdots & & & 0 \\
\tbinom{r}{3} & -\tbinom{r}{2} & r & & & &  0\\
\vdots & & & \ddots  & &  \vdots & \vdots \\
~\\
0 & 0&0 & \cdots \cdots && r & 0 \\ 
 \end{pmatrix} \, .
\eeq

The following is then immediate.

\begin{prop}\label{k1lens}
The $(n+1) \times (n+1)$ matrix $A$ has rank $n$, whence
$$
K_1(C(\rL^{(n,r)}_q)) \simeq \Z \, .
$$
\end{prop}

On the other hand, 
the structure of the cokernel of the matrix $A$  
depends on the divisibility properties of the integer $r$. 
Since $\coker(A) \simeq \Z^{n+1} / \mathrm{Im}(A)$ and $\mathrm{Im}(A)$ being generated by the columns of $A$, the vanishing of these columns yields conditions on the generators making them torsion classes in general. 
Indeed, upon pulling back to the lens space, 
the vanishing of the the $j$-th column is just the condition that the pulled back line bundles satisfy $\tilde{\LL}_{-(r+j)}=\tilde{\LL}_{-j}$; thus this vanishing contains geometric information. 

However, to quickly determine $\coker(A)$ (although not directly its generators) one can
use the {\em Smith normal form} for matrices over a principal ideal domain, such as $\Z$. 
 Thus \cite{Sm} (\textit{cf}. \cite[Thm. 26.2 and Thm. 27.1]{M}) there exist invertible matrices $P$ and $Q$ having integer entries which transform $A$ to a diagonal matrix 
\beq\label{dsf}
\mathrm{Sm}(A) := PA Q = \diag (\alpha_1, \cdots, \alpha_n,0 ) \, ,
\eeq
with integer entries $\alpha_i \geq 1$, ordered in such a way that $\alpha_i \mid \alpha_{i+1}$ for $1 \leq i \leq n$.
These integers are algorithmically and explicitly given by 
$$
\alpha_1= d_1(A) \, , \qquad \alpha_i= {d_i(A)} / {d_{i-1}(A)} \, , \quad \mathrm{for~ each} \quad 2 \leq i \leq n \, ,
$$ 
where $d_i(A)$ is the greatest  common divisor of the non-zero determinants of the minors of order $i$ of the matrix 
$A$. The above leads directly to the following.
 
\begin{prop}\label{k0lens}
It holds that 
$$
\coker(A) \simeq \coker(\mathrm{Sm}(A)) = \Z \oplus \Z / {\alpha_1} \Z \oplus \dots \oplus \Z / \alpha_n\Z \, .
$$
As a consequence,
$$
K_0(C(\rL^{(n,r)}_q)) \simeq \Z \oplus \Z_{\alpha_1} \oplus \dots \oplus \Z_{\alpha_n} \, ,  
$$
with the convention that $\Z_1 = \Z / 1 \Z$ is the trivial group.  
\end{prop}

As already mentioned, the merit of our construction is not in the computation of the K-theory groups -- these are found for instance by using graph algebras 
as in \cite{HS02}. Owing to the explicit diagonalization as in \eqref{dsf} and to Prop.~\ref{ku}, we also obtain explicit generators as integral combinations of powers of the pull-back to the lens space $\rL^{(n,r)}_q$ of the generator $u := 1 - [\LL_{-1}]$.
We show how this works and compute $K_0(C(\rL^{(n,r)}_q))$ in some examples.  

\begin{exa}\label{r=2}
\textup{
If $r=2$ one computes $\alpha_1=\alpha_2 = \cdots =\alpha_{n-1}=1$ and $\alpha_n= 2^n$. 
Hence for $\rL^{(n,2)}_q = \rS^{2n+1}_q / \Z_2 = \R P^{2n+1}_q$, the quantum real projective space, we get 
$$
K_0(C(\R P^{2n+1}_q)) = \Z \oplus \Z_{2^n} \, ,   
$$
in agreement with \cite[\S4.2]{HS02} (with a shift $n\to n+1$ from there to here). 
Moreover, we can construct explicitly the generator of the torsion part of the K-theory group. 
We claim this is given by $1-[\tilde{\mathcal{L}}_{-1}]$. 
First of all, owing to $\tilde{\mathcal{L}}_{-2} \simeq \tilde{\mathcal{L}}_{0}$ one has
$$
(1- [\tilde{\mathcal{L}}_{-1}])^2 = 2(1- [\tilde{\mathcal{L}}_{-1}]),
$$
and iterating:
$$
(1- [\tilde{\mathcal{L}}_{-1}])^k = 2^{k-1} (1- [\tilde{\mathcal{L}}_{-1}]) . 
$$
Thus, in a sense one can switch from multiplicative to additive notation.  
Furthermore, from Prop.~\eqref{ku} we know that $u^{n+1}=0$, with $u=1-[\mathcal{L}_{-1}]$. 
When pulled back to the lens space, owing to $\tilde{\mathcal{L}}_{2N} \simeq \tilde{\mathcal{L}}_{0}$ 
and $\tilde{\mathcal{L}}_{2N+1}\simeq \tilde{ \mathcal{L}}_{-1}$, this implies that 
$$
0 = (1-[\tilde{\mathcal{L}}_{-1}])^{n+1} = 2^{n} (1 - [\tilde{\mathcal{L}}_{-1}]) .
$$
This amounts to saying that the generator $1-[\tilde{\mathcal{L}}_{-1}]$ is cyclic with the correct order $2^{n}$.
}
\end{exa}
 
\begin{exa} 
\textup{
For $n=1$ there is only one $\alpha_1=r$. Then in this case one has
$$ 
K_0 (C(\rL^{(1,r)}_q)) = \Z \oplus \Z_r \, .  
$$ 
From its very definition $[\widetilde{\LL}_{-r}]=1$, thus  $\widetilde{\LL}_{-1}$ generates the torsion part. Alternatively, from $u^2=0$ 
it follows that  $\LL_{-j} = -(j-1) + j \LL_{-1}$ for all $j>0$; upon lifting to $\rL^{(1,r)}_q$, for $j=r$ this yields $r(1-[\widetilde{\LL}_{-1}]) = 0$,
that is $1-[\tilde{\mathcal{L}}_{-1}]$ is cyclic of order $r$.       
}
\end{exa}
\begin{exa} 
\textup{
For $n=2$ there are two cases, according to whether $r$ is even or odd.
For the $\alpha$'s in Prop.~\ref{k0lens} one finds:
$$
(\alpha_1, \alpha_2) = 
\begin{cases} 
(r / 2, 2r) & \qquad \mbox{if} \quad r \quad \mbox{even} \\
(r, r) & \qquad \mbox{if}  \quad r \quad \mbox{odd} 
\end{cases} \, .   
$$
As a consequence one has that
$$ 
K_0 (C(\rL^{(2,r)}_q)) = \begin{cases} \Z \oplus \Z_{\frac{r}{2}} \oplus \Z_{2r} 
& \qquad \mbox{if} \quad r \quad \mbox{even} \\
\Z \oplus \Z_{r} \oplus \Z_{r} 
& \qquad \mbox{if}  \quad r \quad \mbox{odd} 
\end{cases} \, .    
$$
This is in agreement with \cite[Prop.~2.3]{HS03} (once again with a shift $n\to n+1$). In particular, for $r=2$ we get back the case of Example~\ref{r=2}. 
In order to identify generators in the two cases, we start from $[\widetilde{\LL}_{-r}]=1$. Direct computations from the 
conditions $[\widetilde{\LL}_{-(r+j)}] = [\widetilde{\LL}_{-j}]$ for $j=0, \cdots, r-1$ lead to 
\beq\label{2cond}
\tfrac{1}{2} r (r-1) \, \tilde{u}^{2} - r \, \tilde{u} = 0 \quad \textup{and} \quad r \, \tilde{u}^{2}=0 \, ,
\eeq
where $\tilde{u}=1-[\tilde{\mathcal{L}}_{-1}]$. Indeed these are just the lifts to the lens space $\rL^{(2,r)}_q$ of the non-vanishing columns of the corresponding matrix $A$ in \eqref{matalp}. \\~\\
When $r=2k$ is even, we have conditions coming from $(\widetilde{\LL}_{-2})^k\simeq \widetilde{\LL}_{0}$. In fact, due to $[\widetilde{\LL}_{-2k}]=1$, one has
$(1-[\tilde{\mathcal{L}}_{-k}])^2= 2 (1 - [\tilde{\mathcal{L}}_{-k}])$, leading to
$$
0 = (1-[\tilde{\mathcal{L}}_{-k}])^{3} = 4 (1 - [\tilde{\mathcal{L}}_{-k}])= 4k \, \tilde{u} - 2k(k-1) \, \tilde{u}^{2}. 
$$
Together with the conditions \eqref{2cond} this yields 
$$
\tfrac{1}{2} r  \, ( \tilde{u}^{2} + 2 \, \tilde{u} )=0 \quad \textup{and} \quad 2 r \, \tilde{u} =0 , 
$$
that is $\tilde{u}^{2} + 2 \, \tilde{u}$ is of order $r / 2$ while $\tilde{u}$ is of order $2r$ (again, for $r=2$ this is consistent with 
the result of Example~\ref{r=2}, the first `generator' collapsing to the condition $\tilde{u}^{2} + 2 \, \tilde{u}=0$). \\~\\
When $r=2k+1$ is odd, the conditions \eqref{2cond} just say that $\tilde{u}$ and $\tilde{u}^{2}$ are cyclic of order $r$: 
$$
r \, \tilde{u} =0 \quad \textup{and} \quad r \, \tilde{u}^{2}=0 .
$$
}
\end{exa}

\begin{exa}
\textup{
When $n=3$ the selection of generators for the torsion groups is more involved but still `doable'. 
We compute explicitly in App.~\ref{coke} the cokernel of the matrix $A$ in \eqref{matalp} and list here 
the K-theory groups as well as the generators obtained by lifting to the lens space the cokernel of $A$ 
via the surjective map $j_*$. As before we denote $\tilde{u}=1-[\tilde{\mathcal{L}}_{-1}]$. 
There are now four possibilities. For the $\alpha$'s in Prop.~\ref{k0lens} one finds:
\begin{center}
\begin{tabular}{c|c|c|c|c|}
& $ 6 \mid r $ & $2\mid r, 3 \nmid r$ &$ 2\nmid r, 3 \mid r$ & $2 \nmid r, 3 \nmid r$\\
\hline \hline
$ \alpha_1$ & $r/6$& $r/2$ &  $r/3$ &$r$\\
 \hline
 $\alpha_2$ & $r/2$&$r/2$ & $r$ & $r$ \\
\hline  
  $\alpha_3$ & $12r$& $4r$& $3r$ & $r$ \\
\end{tabular} \, .
\end{center}
As a consequence: \\~\\
\noindent
Case $r\equiv 0 \pmod 6$:
$$ 
K_0 (C(\rL^{(3,r)}_q)) = \Z \oplus \Z_{\frac{r}{6}} \oplus \Z_{\frac{r}{2}} \oplus \Z_{12r}
$$
with generators
$$
\tilde{u}^{3} + 12 \, \tilde{u} \, , \quad \tilde{u}^{2} + 6 \, \tilde{u} \, , \quad \tilde{u} \, , 
$$
of order $r / 6$, $ r / 2 $ and ${12r}$, respectively. For the particular case $r=6$, the first torsion part is absent, 
one has $\tilde{u}^{3} + 12 \, \tilde{u} = 0$, and
$$
K_0 (C(\rL^{(3,6)}_q)) = \Z \oplus \Z_{3} \oplus \Z_{72} . 
$$
\\~
\noindent
Case $r\equiv 2,4 \pmod 6$:
$$
K_0 (C(\rL^{(3,r)}_q)) = \Z \oplus \Z_{\frac{r}{2}} \oplus \Z_{\frac{r}{2}} \oplus \Z_{4r}
$$
with generators
$$
\tilde{u}^{3} + 2 \, \tilde{u}^{2} \, , \quad  \tilde{u}^{2} + 2 \, \tilde{u}  \, , \quad \tilde{u} \, ,
$$
of order $r / 2$, $ r / 2 $ and ${4r}$, respectively. The particular case $r=2$ goes back to Example~\ref{r=2} with the first and second torsion parts absent and the condition $\tilde{u}^{2}+2 \, \tilde{u}=0$ as in there. 
\\~\\ 
Case $r\equiv 3 \pmod 6$:
$$
K_0 (C(\rL^{(3,r)}_q)) = \Z \oplus \Z_{\frac{r}{3}} \oplus \Z_{r} \oplus \Z_{3r}
$$
with generators
$$
\tilde{u}^{3} + 3 \, \tilde{u} \, , \quad  \tilde{u}^{2}  \, , \quad \tilde{u} \, ,
$$
of order $r / 3$, $r$ and ${3r}$, respectively. For the particular case $r=3$ the first torsion part is absent, one has 
$\tilde{u}^{3} + 3 \, \tilde{u} = 0$, and
$$
K_0 (C(\rL^{(3,3)}_q)) = \Z \oplus \Z_{3} \oplus \Z_{9} . 
$$
\\~
\noindent
Case $r\equiv 1,5 \pmod 6$:
$$ 
K_0 (C(\rL^{(3,r)}_q)) = \Z \oplus \Z_{r} \oplus \Z_{r} \oplus \Z_{r} 
$$
with the three generators of order $r$ given by
$$
\tilde{u}^{3} \, , \quad  \tilde{u}^{2}  \, , \quad \tilde{u} \, .
$$
}
\end{exa}
To further illustrate the construction, we mention the next case of the dimension $n$, for which we list the K-theory groups.
\begin{exa}
\textup{
When $n=4$ there are $8$ possibilities. For the $\alpha$'s in Prop.~\ref{k0lens} one finds: \\
\begin{center}
\begin{tabular}{c|c|c|c|c|c|cl}
& $ 24 \mid r $ & $12\mid r; 8\nmid r$ & $ 8 \mid r; 6 \nmid r,$ & $6 \mid r; 4\nmid r, $ & $4 \mid r; 3, 8 \nmid r$ \\
\hline \hline
$ \alpha_1$ & $r/24$& $r/12$ &  $r/8$ &$r/6$ &$r/4 $ & \\
 \hline
$ \alpha_2$ & $r/6$& $r/12$ &  $r/4$ &$r/6$ &$r/4 $  & \\
\hline
$ \alpha_3$ & $6r$& $12r$ &  $4r$ &$4r$ &$2r $ & \\
\hline 
$ \alpha_4$ & $24r$& $12r$ &  $8r$ & $12r$ &$8r $ & \\
\hline
\end{tabular} \quad .... \qquad\qquad \qquad
\end{center}
\begin{center}
 \qquad\qquad  \qquad\qquad  \qquad\qquad  \qquad  .... \quad \begin{tabular}{c|c|c|c|}
 & $3\mid r; 2 \nmid r$& $2 \mid r; 3,4 \nmid r$ & $2\nmid r; 3 \nmid r$\\
\hline \hline
 & $r/3$ & $r/2$ & $r$\\
 \hline
 & $r/3$ & $r/2$ & $r$\\
\hline
& $r$ & $r/2$ & $r$\\
\hline 
& $9r$ & $8r$ & $r$ \\
\hline
\end{tabular} \, .
\end{center}
As a consequence,
$$ 
K_0 (C(\rL^{(4,r)}_q)) = \begin{cases} \Z \oplus \Z_{\frac{r}{24}} \oplus \Z_{\frac{r}{6}} 
\oplus \Z_{6r} \oplus \Z_{24r}  & r\equiv 0 \pmod {24} \\ 
\Z \oplus \Z_{\frac{r}{12}} \oplus \Z_{\frac{r}{12}} \oplus \Z_{12r} \oplus \Z_{12r}  & r\equiv 12 \pmod {24} \\ 
\Z \oplus \Z_{\frac{r}{8}} \oplus \Z_{\frac{r}{4}} \oplus \Z_{4r} \oplus \Z_{8r}  & r\equiv 8, 16 \pmod {24} \\ 
\Z \oplus \Z_{\frac{r}{6}} \oplus \Z_{\frac{r}{6}} \oplus \Z_{4r} \oplus \Z_{12r}  & r\equiv 6 \pmod {12} \\ 
\Z \oplus \Z_{\frac{r}{4}} \oplus \Z_{\frac{r}{4}} \oplus \Z_{2r} \oplus \Z_{8r}  & r\equiv 4,20 \pmod {24} \\ 
\Z \oplus \Z_{\frac{r}{3}} \oplus \Z_{\frac{r}{3}} \oplus \Z_{r} \oplus \Z_{9r}  & r\equiv 3, 9 \pmod {12} \\ 
\Z \oplus \Z_{\frac{r}{2}} \oplus \Z_{\frac{r}{2}} \oplus \Z_{\frac{r}{2}} \oplus \Z_{8r}  & r\equiv 2 \pmod {12} \\ 
\Z \oplus \Z_{r} \oplus \Z_{r} \oplus \Z_{r} \oplus \Z_{r}  & r\equiv 1,5 \pmod {6} \\ 
\end{cases} \, .  
$$
}
\end{exa}
~\\
\noindent
We leave to the diligent reader the determination of the corresponding generators.

\section{Final remarks and future perspectives}

In the classical geometry of `commutative' spaces, the Gysin sequences for principal $\U(1)$-bundles (at least at the level of de Rham cohomology) play a particularly important role both in T-duality  theories and in Chern-Simons theory; these are topics of intensely active current research, although they by no means exhaust the uses of such a sequence. 

In T-duality theories one studies principal $\U(1)$-bundles $E\to M$ over a compact manifold $M$ and their associated line bundles. The Gysin sequence relates the $H$-flux (a given three-form on the total space $E$) to the curvature (a two-form on the base space $M$) of a connection on a dual line bundle $E'$. The curvature is indeed though of as representing the Chern class of the dual bundle $E'$ on $M$. In the case of a two-dimensional base manifold $M$ this also gives an isomorphism between Dixmier-Douady classes on $E$ and line bundles on $M$ (\textit{cf}. pages 385 and 391 of the seminal paper \cite{T-D}).  

In Chern-Simons theory, the importance of the role of the Gysin sequence is in the evaluation of the path integral on $\U(1)$-bundles over smooth curves, where it facilitates the counting of those  $\U(1)$-bundles over the total space which arise as pull-backs from the base (\textit{cf}. page 26 of the fundamental paper \cite{C-S}). Thus in both of these applications of the Gysin sequence, what is important is the use of explicit representatives of elements in cohomology, rather than the simple knowledge of the cohomology groups themselves.

In the present paper we have presented a novel geometric approach to quantum lens spaces and their noncommutative topology. 
This is done via a Gysin sequence in K-theory for these spaces which, to the best of our knowledge, have never been studied before in this context.
The strength of our construction is not only the matter of computing the K-theory groups, which could and has be done by means of graph algebras, although ours is a method of a novel sort. Our central result is that we also obtain explicit geometric generators as classes of `line bundles' which generically are torsion ones.  
These line bundles are pulled-back from line bundles over the quantum projective spaces, the `base spaces' of $\U(1)$-bundles whose `total spaces' are the quantum lens spaces themselves. 

There are several potential but  important 
applications to problems which could not be treated with the old methods, that is to say, where explicit geometric generators for (classes of) bundles play a crucial role. For our scientific taste, applications concern mainly $T$-duality theory and Chern-Simons theory for quantum spaces: they are under current investigation and will be reported elsewhere. 

\addtocontents{toc}{\protect\setcounter{tocdepth}{1}}

\appendix
 
\section{Principal bundle structures}\label{ncbundles}
In this Appendix we carry out the proofs of Prop.~\ref{pr-u1t} and Prop.~\ref{pr-Zr}. While for the principal bundles there the structure groups are ordinary (Abelian) groups, $\Z_r$ and $\widetilde{\U}(1)$ respectively, we use their dual Hopf algebras that give direct and easy algebraic methods. 

A {\em noncommutative principal bundle} is a triple $(\A,H,\mathcal{F})$, where $\A$ is the $*$-algebra of functions on `the total space',  $H$ is the Hopf $*$-algebra of functions on the `structure group', with $\A$ being a right (say) $H$-comodule $*$-algebra, that is there is a right coaction
$$
\Delta_R:\A\to \A\otimes H . 
$$
The functions on `the base space' are given by the $*$-subalgebra of coinvariant elements:
$$
\mathcal{F}:=\{a\in \A~|~\Delta_R(a)=a\otimes 1\} .
$$
Conditions to be satisfied are imposed via a suitable sequence.  
Let $\Omega^1_{un}(\B)$ denote the bimodule of universal differential forms over a unital algebra $\B$ and $\epsilon_{H}$ be 
the counit of the Hopf algebra $H$. Principality of the bundle is expressed by requiring the sequence 
\beq
\label{exactseqbndl}
0 \to \mathcal{A}(\Omega^1_{un}(\mathcal{F})) \mathcal{A} \to \Omega^1_{un}(\mathcal{A}) \xrightarrow{\textup{ver}} \mathcal{A} \otimes \ker \epsilon_{H} \to 0 
\eeq
to be exact. Here the first  map is inclusion while the second one, $\textup{ver}(a\otimes b):=(a\otimes 1)\Delta_R(b)$, generates `vertical one-forms'.
When $H$ is cosemisimple and has an invertible antipode, exactness of the sequence \eqref{exactseqbndl} is equivalent to the statement that the canonical map 
\beq
\label{chi}
\chi: \A\otimes_\mathcal{F}\A \rightarrow \A \otimes H,\qquad 
\chi(a\otimes b):=(a\otimes 1)\Delta_R(b),
\eeq
is an isomorphism (sometimes this is also known as the statement that the triple $(\A,H,\mathcal{F})$ is a {\em Hopf-Galois extension}). Furthermore, things are easier for a cosemisimple Hopf algebra $H$ with bijective antipode, since then the map \eqref{chi} is injective whenever it is surjective and thus it is enough to check surjectivity \cite[Thm.~I]{Sch}.

\subsection{Proof of Prop.~\ref{pr-Zr}} 
The algebra $\A(\Z_r)$ of functions on the cyclic group $\Z_r$ is the $*$-algebra generated by a single element $\zeta$ modulo the relation $\zeta^r=1$. The Hopf  structures are given by coproduct $\Delta(\zeta)=\zeta\otimes\zeta$, counit $\epsilon(\zeta)=1$ and antipode $S(\zeta)=\zeta^*$.

Thus the algebra $\A (\rL^{(n,r)}_q)$ can be also obtained as the algebra of coinvariant elements with respect to the coaction of $\A (\Z_r)$ on $\A(\rS^{2n+1}_q)$ defined on generators by
$$
\Delta_R (z_i) = z_i \otimes \zeta, \qquad \Delta_R (z^*_i) = z^*_i \otimes \zeta^{*}
$$
and extended as algebra map to the whole of $\A(\rS^{2n+1}_q)$. 

The Hopf algebra $\A(\Z_r)$ is certainly cosemisimple and so, to show that the datum 
$(\mathcal{A}(\rS^{2n+1}_q), \A (\Z_r), \A(\rL^{(n,r)}_q))$ is a principal bundle, it is enough to establish surjectivity of the map $\chi$ defined as in \eqref{chi}. For this we use a strategy borrowed from \cite{LPR06}.

Writing $\A=\A(\rS^{2n+1}_q)$, $H=\A(\Z_r)$ and $\mathcal{F}=\A(\rL^{(n,r)}_q)$, a generic element in $\A \otimes H$ is a sum of elements of the form $f \otimes \zeta^{*N}$ with $N =0, \dots, r-1$ and $f \in \A(\rS^{2n+1}_q)$. 
By left $\A(\rS^{2n+1}_q)$-linearity of $\chi$, it is enough to exhibit a pre-image for elements of the form $1 \otimes \zeta^{*N}$, since
if $\gamma \in \A\otimes_\mathcal{F}\A$ is such that $\chi(\gamma) = 1 \otimes \zeta^{*N}$, then $\chi (f\gamma) =f (1 \otimes \zeta^{*N}) = f \otimes \zeta^{*N}$. 

Let $\psi^N_{j_0, \dots, j_n}$ be the vector-valued functions given in \eqref{psi}, and define the element  
\beq\label{pregam}
\gamma := \sum_{j} \psi^{N*}_{j} \otimes \psi^{N}_j ,
\eeq
which is clearly in $\A \otimes_{\mathcal{F}} \A$. 
Denote $\beta^N_j := [j_0,\dots, j_n]! q^{-\sum_{r<s}j_r j_s}$, to lighten notation. Upon applying $\chi$ one obtains
\begin{align*}
\chi(\gamma) &=  \sum_j \chi (\psi^{N*}_{j} \otimes \psi^{N}_{j} ) 
= \sum_{j_0 + \dots + j_n = N} \beta^{N}_{j} ( z_n^{j_n} \dots z_0^{j_0} \otimes 1) 
\, \Delta_R \left( (z_0^{j_0})^* \dots (z_n^{j_n})^* \right) \\ 
&= \sum_j  \beta^{N}_{j} \left( z_{n}^{j_n} \dots z_{0}^{j_0} \otimes 1 \right) 
\cdot \left(  (z_0^{j_0})^{*} \dots (z_n^{j_n})^{*} \otimes \xi^{*N}\right) \\
&=  \sum_j  \beta^N_j \left( z_n^{j_n} \dots z_0^{j_0}  (z_0^{j_0})^* \dots (z_n^{j_n})^*   \otimes \xi^{*N}\right) \\
&=  \left( \sum_j  \beta^N_j  z_n^{j_n} \dots z_0^{j_0}  (z_0^{j_0})^* \dots (z_n^{j_n})^*  \right) \otimes \xi^{*N}  \\
&= 1 \otimes \xi^{*N},
\end{align*}
which is all one needs for proving surjectivity of the map $\chi$.

\subsection{Proof of Prop.~\ref{pr-u1t}}
Let  $\A(\U(1) = \C[\xi, \xi^*] / \langle \xi^*\xi -1 \rangle$ denote the coordinate algebra of the group $\U(1)$.
With $\widetilde{\U}(1) :=\U(1)/\Z_r$, the corresponding coordinate algebra, 
$$
 \A (\widetilde{\U}(1)) := \A(\U(1)/\Z_r) = \A(\U(1))^{\Z_r} ,
$$
is the Hopf $*$-subalgebra of $\A(\U(1))$ generated by the powers $\xi^r$ and $\xi^{*r}$.

Denote $\A'= \A (\rL^{(n,r)}_q)$, $H'=\A(\widetilde{\U}(1))$ and $\mathcal{F}' = \A(\CP^{n}_q)$. As before, 
the datum $(\A', H', \mathcal{F}')$ is seen as a quantum principal bundle via surjectivity of the canonical map 
$$
\chi' : \A' \otimes_{\mathcal{F}'} \A' \rightarrow \A' \otimes H', 
$$
surjectivity proved again by exhibiting a pre-image for elements of the kind $1 \otimes \xi^{*rN}$.

For this, we observe that the vectors $\psi^N_{j_0, \dots, j_n}$ in \ref{psi} have entries in $\A' = \A(\rL^{(n,r)}_q)$ precisely when $N$ is a multiple of $r$. Then, in parallel with \eqref{pregam}, the element of $\A' \otimes_{\mathcal{F}'} \A'$,
\beq\label{pregamr}
\gamma' := \sum_{j} \psi^{rN*}_{j} \otimes \psi^{rN}_j ,
\eeq
is mapped to $1 \otimes \zeta^{*rN}$ by the canonical map $\chi'$; enough for the surjectivity of the latter.

\section{Computing cokernels}\label{coke}

We compute explicitly the cokernel of the matrix $A$ in \eqref{matalp} when $n=3$. 
This is a bit involved and, depending on the divisibility properties of the integer $r$, requires considering different cases for $r$. 
Now $\coker(A) \simeq \Z^{4} / \mathrm{Im}(A)$; since $\mathrm{Im}(A)$ is generated by the columns of $A$, 
the vanishing of these columns yields conditions on the generators of $\coker(A)$.

\medskip
\noindent $r=6k$. \quad
The columns of the matrix $A$ yield the constraints
$$
\begin{cases}
6k u - 3k(6k-1) u^{2} + k(6k-1)(6k-2) u^{3}  =  0\, ,\\
6k u^{2} - 3k(6k-1)u^{3}  = 0\, ,\\
6k u^{3}  = 0 \, .
\end{cases}
$$
Substituting the third equation into the first and second ones yields the relations
$$
\begin{cases}
6k u + 3k ( 1 - 6k ) u^{2} + 2 k u^{3}  =  0\, ,\\
3k (2u^{2} \pm u^{3}) = 0 \quad \Rightarrow \quad 12 k u^{2} = 0\, , \\
6k u^{3}  = 0 \, .
\end{cases}
$$
By multiplying the first equation by two and using the second, one gets
$$
k(12 u+u^{3}) = 0 ,
$$
that is $12 u+u^{3}$ has order $k=r/6$.  On the other and, the first equation can be multiplied by three,
after which the use of the third equation yields $18 k u + 9k ( 1 - 6k) u^{2} = 0$. 
Now, modulo $12 k u^2$, one has $9k ( 1 - 6k ) u^{2} \equiv 3 k u^{2}$, which transforms the previous equation into
$$
3k (6u + u^{2}) = 0, 
$$
that is $(6u + u^{2})$ has order $3k=r/2$. Finally, multiplying the first equation by $6$ and using $6k u^{3}  = 0$ or the second equation by $4$ and using $12 k u^{2} = 0$, it follows that
$$
72 k u=0,
$$
i.e. $u$ has order $12 r$.

\medskip
\noindent $r=6k+2$ and $r=6k-2$. \quad
For the first case, the columns of $A$ yield the constraints
$$
\begin{cases}
2 (3k+1) u-(3k+1)(6k+1)u^{2}+2k(3k+1)(6k+1)u^{3} =  0 \, ,\\
2 (3k+1) u^{2}-(3k+1)(6k+1)u^{3}  =  0 \, ,\\
2 (3k+1) u^{3}  = 0 \, .
\end{cases}
$$
These can be rewritten as
$$
\begin{cases}
2 (3k+1) u-(3k+1)(6k+1)u^{2} =  0\, ,\\
2 (3k+1) u^{2} \mp (3k+1)u^{3}  = 0 \quad \Rightarrow \quad 4 (3 k +1) u^{2} = 0\, , \\
2 (3k+1) u^{3}  =  0
\end{cases}
$$
and from the second equation one immediately gets
$$
(3k+1)(2u^2+u^{3})=0,
$$
which says that $(2u^{2}+u^{3})$ has order $3k+1 = r/2$. On the other hand, 
modulo $4 (3 k +1) u^{2}$ one has that $(3k+1)(6k+1) u^{2} \equiv (-1)^k (3k+1) u^{2}$ which transforms the first equation to 
$$
(3k+1)(2u + (-1)^{k+1} u^{2})=0,
$$
that is $2u + (-1)^{k+1} u^{2}$ has order $3k+1 = r/2$. Moreover, using again  $4 (3 k +1) u^{2}=2$ this also yields
$$
4(3k+1) 2 u = 0, 
$$
that is $u$ has order $8(3k+1) = 4 r$.  

Analogous computations and results holds for $r=6k-2$. 

\medskip
\noindent $r=6k+3$. \quad
The columns of the matrix $A$ yield the constraints$$
\begin{cases}
3 (2k+1) u- 3(3k+1)(2k+1) u^{2}+(2k+1)(3k+1)(6k+1)u^{3} = 0\, ,\\
3 (2k+1) u^{2}-3(3k+1)(2k+1) u^{3}  = 0\, ,\\
3 (2k+1) u^{3} = 0 \, .
\end{cases}
$$
These can be rewritten as
$$
\begin{cases}
3 (2k+1) u+(2k+1)(3k+1)(6k+1)u^{3} =  0\, ,\\
3 (2k+1) u^{2} = 0\, ,\\
3 (2k+1) u^{3} = 0\, ,
\end{cases}
$$
that is both $u^{3}$ and $u^{2}$ have order $3 (2k+1) = r$. Moreover, using twice the last equation in the first one leads to
\begin{eqnarray*}
0 & = & (6k+3)u-2(2k+1)(3k+1)u^{3}=(2k+1)(3u-(6k+2)u^{3}) \\
 & = & (2k+1)(3u+u^{3}),
\end{eqnarray*}
which says that $3u+u^{3}$has order $2k+1=r/3$. Finally 
$$
9(2k+1) u=3(2k+1)u^{3}=0,
$$
hence $u$ has order $9(2k+1)=3r$.

\medskip
\noindent $r=6k+1$ and $r=6k-1$. \quad
The columns of the matrix $A$ yield the constraints
$$
\begin{cases}
(6k+1)u-3k(6k+1)u^{2}+k(6k+1)(6k-1)u^{3} =  0\, ,\\
(6k+1)u^{2}-3k(6k+1)u^{3} = 0\, ,\\
(6k+1)u^{3} = 0 \, .
\end{cases}
$$
These just tell us that $u,u^{2},u^{3}$ all have order $6k+1=r$. 

Analogous computations and results hold for $r=6k-1$.


\begin{thebibliography}{9999}

\bibitem{BJ83}
S.~Baaj, P.~Julg, \emph{Th\'eorie Bivariante de Kasparov et Operateurs non-Born\'{e}s dans les C*-Modules Hilbertiens}, C.R. Acad. Sc. Paris 296 (1983), 875--878.

\bibitem{Bl}
B.~Blackadar, \textit{$K$-Theory for Operator Algebras}, 2nd edition,
Cambridge University Press, 1998.

\bibitem{C-S}
M.~Blau, G.~Thompson,
\emph{Chern-Simons theory on $S^1$-bundles: abelianisation and q-deformed Yang-Mills theory},
JHEP 05 (2006), 003.

\bibitem{T-D}
P.~Bouwknegt, J.~Evslin, V.~Mathai, 
\emph{T-Duality: Topology Change from $H$-Flux}
Commun. Math. Phys. 249 (2004), 383--415. 

\bibitem{BF12}
T.~Brzezinski, S.A. Fairfax, \emph{Quantum Teardrops}, Commun. Math. Phys. 316 (2012), 151--170.

\bibitem{BM}
T.~Brzezinski, S.~Majid, \emph{Quantum Group Gauge Theory on Quantum Spaces}, Commun. Math. Phys. 157
(1993), 591--638; Erratum 167 (1995), 235.

\bibitem{CPR}
A.L.~Carey, J.~Phillips, A.~Rennie, \textit{Noncommutative Atiyah-Patodi-Singer Boundary Conditions and Index Pairings in KK-Theory}, J. reine angew. Math. 643 (2010), 59--109.

\bibitem{CNNR08}
A.L.~Carey, S.~Neshveyev, R.~Nest, A.~Rennie, \textit{Twisted Cyclic Theory, Equivariant KK-theory and KMS States}, J. reine angew. Math. 650 (2011), 161--191.

\bibitem{Cu84}
J.~Cuntz, \textit{On the Homotopy Groups for the Space of Endomorphisms of a $C^*$-Algebra},
in: ``Operator Algebras and Group Representations", Pitman, London, 1984, 124--137.

\bibitem{DD09}
F.~D'Andrea, L.~D\c{a}browski, \textit{Dirac Operators on Quantum Projective Spaces},  Commun. Math. Phys. 295 (2010), 731--790.

\bibitem{DL09}
F.~D'Andrea, G.~Landi, \textit{Anti-Self-Dual Connections on the Quantum
 Projective Plane: Monopoles}, Commun. Math. Phys. 297 (2010), 841--893.

\bibitem{DL10}
\bysame,
\textit{Bounded and Unbounded Fredholm Modules for Quantum Projective Spaces},
J. K-Theory 6 (2010), 231--240.

\bibitem{DL13}
\bysame, \textit{Geometry of Quantum Projective Spaces},
in: ``Noncommutative Geometry and Physics 3". G. Dito et al. Editors, 
Keio COE Lecture Series on Mathematical Science: Vol. 1. 
World Scientific, Singapore, 2013, pp. 373--416.

\bibitem{HRZ11}
P.M.~Hajac, A.~Rennie, B.~Zieli{\'n}ski, \emph{The K-Theory of Heegaard Quantum Lens Spaces}, J. Noncommut. Geom. 7 (2013), 1185--1216.

\bibitem{HL04}
E.~Hawkins, G.~Landi, \textit{Fredholm Modules for Quantum Euclidean
 Spheres}, J.~Geom. Phys. 49 (2004), 272--293.

\bibitem{HR}
N.~Higson, J.~Roe, \textit{Analytic K-Homology}. Oxford University Press, 2000

\bibitem{HS02}
J.H.~Hong, W.~Szyma{\'n}ski, \textit{Quantum Spheres and Projective Spaces
 as Graph Algebras}, Commun. Math. Phys. 232 (2002), 157--188.

\bibitem{HS03}
\bysame, \textit{Quantum Lens Spaces and Graph Algebras}, Pac. J. Math. 211 (2003), 249--263.

\bibitem{Ka78} M.~Karoubi, \textit{K-Theory: an Introduction}. Grundlehren der math. Wiss.  226, Springer, 1978.

\bibitem{Kas} 
G.G.~Kasparov, The Operator $K$-Functor and Extensions of
$C^{*}$-Algebras, {\em Izv. Akad. Nauk. SSSR Ser. Mat}. 44(1980), 571--636; English
translation, {\em Math. USSR-Izv}. 16 (1981), 513--572.

\bibitem{KLS}
M.~Khalkhali, G.~Landi, W.D.~van Suijlekom, \textit{ Holomorphic Structures on the Quantum Projective Line}, Int. Math. Res. Notices 4 (2011), 851--884.

\bibitem{Lance}
E.C.~Lance, {\em Hilbert C*-Modules: a Toolkit for Operator Algebraists}, London Mathematical Society Lecture Notes Series 210, Cambridge University Press, 1995.

\bibitem{La98}
E.C. Lance, {\em The Compact Quantum Group $SO(3)_q$}, J. Operator Theory, 40 (1998), 295--307.

\bibitem{LPR06}
G.~Landi, C.~Reina, C.~Pagani,
\textit{ A Hopf Bundle over a Quantum Four--Sphere from the Symplectic Group},
Commun. Math. Phys. 263 (2006) 65--88. 

\bibitem{M}
C.C.~Mac Duffee, \textit{The Theory of Matrices}, Courier Dover Publications, 2004.

\bibitem{MT}
K. Matsumoto, J. Tomiyama, 
\textit{Noncommutative Lens Spaces}, 
J. Math. Soc. Japan, 44 (1992), 13--41.

\bibitem{PR}
D.~Pask, A.~Rennie, \textit{The Noncommutative Geometry of Graph $C^*$-Algebras I: the Index Theorem}, J. Func. Anal. 233 (2006), 92--134.

\bibitem{Po95} P. Podle\'s, \textit{Symmetries of Quantum Spaces. Subgroups and Quotient Spaces of Quantum $SU(2)$ and $SO(3)$ Groups}, Commun. Math. Phys. 170 (1995), 1--20.

\bibitem{IP}
I.~Putnam, \textit{An Excision Theorem for the K-Theory of $C^*$-Algebras}, J. Oper. Theory 38 (1997), 151--171.

\bibitem{Sch} 
H.~Schneider, \textit{Principal Homogeneous Spaces for Arbitrary Hopf Algebras}, Israel J. Math. 72 (1990),
167--195.

\bibitem{S-V}
A.~Sitarz, J.J.~Venselaar,
\textit{Real Spectral Triples on 3-dimensional Noncommutative Lens Spaces},
 [arXiv:1312.5690].

\bibitem{Sm}
H.J.S.~Smith, \textit{On Systems of Linear Indeterminate Equations and Congruences}, Phil. Trans. R. Soc. Lond. 151 (1861), 293--326.

\bibitem{VS91}
L.~Vaksman, Ya. Soibelman, \textit{The Algebra of Functions on the Quantum
 Group $\textup{SU}(n+1)$ and Odd-Dimensional Quantum Spheres}, Leningrad Math. J. 2 (1991), 
 1023--1042.

\bibitem{We00}
M.~Welk, 
\textit{Differential Calculus on Quantum Projective Spaces}, 
in: ``Quantum Groups and Integrable Systems (Prague, 2000)", Czech. J. Phys. 50 (2000), 219--224.

\bibitem{Wor88}
S.L.~Woronowicz, \textit{Tannaka-Krein Duality for Compact Matrix Pseudogroups. Twisted $\textup{SU}(N)$ group}, Inv. Math. 93 (1988), 35-76. 

\end{thebibliography}
\end{document}